# On-Orbit Servicing Optimization Framework with High- and Low-Thrust Propulsion Tradeoff

Tristan Sarton du Jonchay[*], Hao Chen[*], Masafumi Isaji[*], Yuri Shimane[*]

and Koki Ho[†,‡]

*Georgia Institute of Technology, Atlanta, GA, 30332, USA*

**This paper proposes an on-orbit servicing logistics optimization framework capable of performing the short-term operational scheduling and long-term strategic planning of sustainable servicing infrastructures that involve high-thrust, low-thrust, and/or multimodal servicers supported by orbital depots. The proposed framework generalizes the state-of-the-art on-orbit servicing logistics optimization method by incorporating user-defined trajectory models and optimizing the logistics operations with the propulsion technology and trajectory tradeoff in consideration. Mixed-Integer Linear Programming is leveraged to find the optimal operations of the servicers over a given period, while the Rolling Horizon approach is used to consider a long time horizon accounting for the uncertainties in service demand. Several analyses are carried out to demonstrate the value of the proposed framework in automatically trading off the high- and low-thrust propulsion systems for both short-term operational scheduling and long-term strategic planning of on-orbit servicing infrastructures.**

## Nomenclature

| | | |
|---|---|---|
| $B_{vst}$ | = | Servicer dispatch variables |
| $\mathcal{C}$ | = | Index set of Customer Nodes |
| $c_s^{\text{delay}}$ | = | Penalty fees per unit time for providing a service with delays |

[*] Ph.D. Student, Daniel Guggenheim School of Aerospace Engineering, Atlanta, GA, AIAA Student Member.
[†] Assistant Professor, Daniel Guggenheim School of Aerospace Engineering, Atlanta, GA, AIAA Senior Member.
[‡] Corresponding author. E-mail address: kokiho@gatech.edu.

$c_v^{\text{dep}}$ = Cost per unit time for operating an orbital depot

$c^{\text{launch}}$ = Launch cost per unit mass

$c_k^{\text{pdm}}$ = Purchase, development, and manufacturing costs of the non-vehicle commodities

$c_v^{\text{pdm}}$ = Purchase, development, and manufacturing costs of the vehicle commodity

$\boldsymbol{d}_{vs}$ = Demand parameters for the non-vehicle commodities at the Customer Nodes

$\mathcal{E}$ = Index set of Earth Nodes

$\boldsymbol{e}_v$ = Vehicle design parameters

$F$ = Thrust force

$g_0$ = Earth's gravitational acceleration at sea level

$H_{vst}$ = Service assignment variables

$i, j$ = Node indices

$I_{sp}$ = Specific impulse

$J_{\text{delay}}$ = Component of the objective function capturing the penalty fees due to service delays

$J_{\text{dep}}$ = Component of the objective function capturing the operating costs of the depots

$J_{\text{launch}}$ = Component of the objective function capturing the launch costs

$J_{\text{pdm}}$ = Component of the objective function capturing the purchase, development, and manufacturing costs

$J_{\text{revenues}}$ = Component of the objective function capturing the revenues generated by the services

$J_{\text{serv}}$ = Component of the objective function capturing the operating costs of the servicers

$\mathcal{K}$ = Index set of commodities

$\mathcal{K}_{\text{cont}}$ = Index set of continuous commodities

$\mathcal{K}_{\text{int}}$ = Index set of integer commodities

$\mathcal{K}_{\text{tools}}$ = Index set of servicer tools

$k$ = Non-vehicle commodity index

$\mathcal{L}_i$ = Union of the service window sets for a particular customer satellite

$\overline{\overline{M}}_{vij}$ = Concurrency matrix

$M_{vijqr}^{\text{ub}}$ = Servicer mass upper bound

$m_k$ = Unit mass for the non-vehicle integer commodities

$m_v$ = Unit mass for the vehicle commodity

$m_{\mathrm{p}}$ = Propellant mass consumption

$m_0$ = Servicer's initial mass

$\mathcal{N}$ = Index set for all nodes in the static network

$\mathcal{N}_{\mathrm{orb}}$ = Index set for the Orbital Nodes in the static network

$\mathcal{P}$ = Index set for the on-orbit servicing parking nodes in the static network

$P^{\pm}_{vijqrt}$ = Commodity flow variables representing the propellant of the servicers

$\mathcal{Q}_v$ = Index set for time of flight options

$\bar{\bar{Q}}'_{vi}$ = Commodity transformation matrix over holdover arcs

$\bar{\bar{Q}}''_{vijqr}$ = Commodity transformation matrix over transportation arcs

$q$ = Time of flight index

$\mathcal{R}_{vq}$ = Index set for trajectory options

$r$ = Trajectory index

$r_s$ = Revenues received upon providing services

$\mathcal{S}$ = Index set for all service needs

$\mathcal{S}_i$ = Index set for the service needs arising at a specific Customer Node

$s$ = Service need index

$\mathcal{T}$ = Index set for the time steps

$T$ = Service time interval in time-expanded network

$t, \tau$ = Time indices

$U^{\pm}_{vijqrtk}$ = Commodity flow variables for the non-vehicle commodities over the transportation arcs

$\mathcal{V}$ = Index set for the vehicles

$\mathcal{V}_{\mathrm{dep}}$ = Index set for the orbital depots

$\mathcal{V}_{\mathrm{serv}}$ = Index set for the servicers

$v$ = Vehicle index

$\mathcal{W}_s$ = Service window index set

$W^{\pm}_{vijqrt}$ = Commodity flow variables for the vehicle commodity over the transportation arcs

$X_{vitk}^{\pm}$ = Commodity flow variables for the non-vehicle commodities over the holdover arcs

$Y_{vit}^{\pm}$ = Commodity flow variables for the vehicle commodity over the holdover arcs

$Z_{vijqrt}^{+}$ = Total mass of the servicers over the transportation arcs

$\beta_{s\tau t}$ = Binary parameters relating servicer dispatch variables and service assignment variables

$\gamma_{sk}$ = Service-tool mapping between service *s* and tool *k*

$\Delta_t$ = Length of holdover arc starting at time step $t$

$\Delta_t'$ = Length of holdover arc ending at time step $t$

$\Delta V$ = Characteristic velocity

$\xi_{ts}$ = Binary parameters relating the time steps to the sets of service windows

$\boldsymbol{\sigma}_{it}$ = Supply of the non-vehicle commodities at the Earth Nodes

$\sigma_{ivt}'$ = Supply of the vehicle commodity at the Earth Nodes

$\tau_s$ = Service need occurrence date

## I. Introduction

Since their inception by science fiction writer Arthur C. Clarke, geosynchronous-equatorial-orbit (GEO) satellites have significantly improved humans' lives on Earth. From telecommunications to weather monitoring, they represent an undeniably critical infrastructure supporting a multitude of terrestrial markets. This, however, comes at the cost of large capital expenditure – on the order of $150 million to $500 million [1] – to manufacture, insure, and launch these large spacecraft to the remote GEO. Until very recently, the traditional paradigm to maintain and upgrade this infrastructure consisted of replacing the outdated satellites with new assets designed to last 15 to 20 years. Besides the obvious issue that satellites designed to operate for that long cannot benefit from rapidly improving technologies, more profound issues exist related to the aversion to risk and innovation. Given the colossal amount of capital needed to deploy and operate space assets, this mindset is understandable but unsustainable if the GEO infrastructure is to become a major player in the cislunar space economy.

On-Orbit Servicing (OOS) may be the tool needed to shift this mindset. OOS has been around for three decades starting with the maintenance of the Hubble Space Telescope by the crews of the Space Shuttle. With the progress of space technologies, U.S. agencies such as NASA and DARPA quickly invested in ambitious robotic servicing missions (*e.g.*, NASA's OSAM 1 [2]; DARPA's Robotic Servicing of Geosynchronous Satellites [3]). Nowadays, as the cost of access to space dramatically decreases, incumbents and new entrants alike are developing innovative technologies and business models that will provide new services to the GEO infrastructure (*e.g.*, refueling, upgrading, repositioning, etc.). With servicing spacecraft (aka *servicers*) regularly visiting GEO assets, satellite manufacturers and operators may change the way they traditionally do business, leading to enhanced resilience to competition and market fluctuations. Northrop Grumman is investing in this new market with its Mission Extension Vehicles (MEV), two of which are already providing life extension services to two of Intelsat's GEO satellites [4]. This is an opportune time for OOS proponents as technology and demand are aligning towards the success of this promising market, valued at $4.5 billion in cumulative revenues by 2028 [5].

Developing a sustainable OOS infrastructure capable of capturing shares of this market will not be without risk. Large initial investments, uncertain demand, and yet-to-be-proven servicing technologies compound that risk. Meticulous strategic planning will be necessary to select the right technologies and business models that ensure a successful OOS endeavor. A typical tradeoff that decision-makers will likely face is between using high-thrust or low-thrust technologies to propel their servicers. Currently, the industry is leaning towards low-thrust propulsion. Northrop Grumman's MEVs for instance integrate Xenon Hall Thrusters [6]. Despite the high price of Xenon gas ($850/kg-$5000/kg in the past 15 years [7,8]), this is a reasonable design given the low service demand and the efficiency and reliability of these thrusters. However, as orbital propellant depots [9] and their routine resupply by cheaper launch vehicles become a reality, and the increase in demand requires more responsive servicing operations, one will naturally question this decision and consider high-thrust propulsion, supported by regular servicer refueling, as a potential alternative. For example, Monomethyl Hydrazine (MMH) is an easily storable liquid propellant whose relatively low price ($170/kg in 2006 [10]) could make high-thrust servicers a competitive option.

To make valuable long-term design decisions, this emerging industry needs tools to explore the OOS design tradespace at the system level and guide technology portfolio management and roadmapping. Such tools will ensure the long-term resilience and robustness of the OOS businesses and their operations despite the uncertainties related to demand and competition. The past literature has already contributed to important aspects of OOS mission planning. Previous projects specifically focused on the rigorous and accurate analysis and design of high-thrust [11] and low-thrust [12] servicer trajectories using high-fidelity force models. Others analyzed the operations of OOS infrastructures from a scheduling [13] and/or design standpoint [14-16] through simulations [15, 17-20], optimization [13, 16, 21], or a mixture of both [14]. None of the above, however, could model and simulate at the system level complex and sustainable OOS infrastructures involving as many orbital depots and servicers as desired over long time horizons and under service demand uncertainties. To tackle that challenge, Ref. [22] proposed to model the GEO orbit as a network and the operations of OOS infrastructures as a time-expanded generalized multi-commodity network problem solved through Mixed-Integer Linear Programming (MILP). In addition, they leveraged the Rolling Horizon (RH) decision making procedure to simulate the operations of OOS infrastructures over the long term and under service demand uncertainties. This state-of-the-art method, however, only considers high-thrust servicers as part of the system-level design tradespace. The modeling of low-thrust servicers was put aside due to unique challenges in integrating non-linear low-thrust trajectory models in a MILP. Nonetheless, with the nascent OOS industry naturally leaning towards using low-thrust engines to propel their servicers, this technology must be accounted for in tradespace analyses.

In response to that background, this paper aims to generalize the state-of-the-art OOS logistics framework proposed in Ref. [22] by modeling and optimizing the operations of not only high-thrust servicers, but also low-thrust and multimodal servicers. A multimodal servicer is modeled as a spacecraft equipped with both high- and low-thrust propulsion systems to find the right balance between service responsiveness and operating costs. The goal of the generalized formulation is to have the optimizer automatically select the best propulsion technology and trajectory of the servicers from a set of discrete options specified by the framework users. We believe this work will give OOS planners a unique insight into one of the most fundamental technological tradeoffs when it comes to in-space transportation: using high-thrust engines, low-thrust engines, or a combination of both to propel the servicers.

Previous works give important methodological pieces to take on the challenge of automated propulsion system tradeoff for the specific problem of OOS planning. This paper uniquely combines them to generalize the state-of-the-art OOS logistics framework. First, Ref. [23] explores the automated tradeoff between high-thrust and low-thrust vehicles for the design of human space exploration campaigns. This was the first time that low-thrust technologies were modeled within space logistics analyses that leverage MILP as the globally optimal solution method. In Ref. [23], the piece-wise linear approximation of the low-thrust model is used along with a new variant of dynamic networks, event-driven networks. In our problem, we do not need the event-driven networks because, unlike the cislunar transportation problem addressed in Ref. [23], the trajectory considered in this paper is for Earth-orbiting satellites to change their phases in the same orbit (i.e. GEO), which can be modeled as a reasonable number of predefined options using time-expanded networks. However, we do need the piecewise linear approximation of non-linear low-thrust propellant consumption models to enable their integration into a MILP. This idea of leveraging piecewise linear approximation of non-linear models in MILP was first proposed by Vielma et al. [24]; it was first introduced to space logistics formulations by Chen and Ho for the integration of non-linear infrastructure design into MILP [25], before Jagannatha and Ho used it to allow the optimizer to automatically perform the tradeoff between high- and low-thrust spacecraft in multi-mission space campaigns [23].

Another methodological piece leveraged in this paper to properly embed the high-thrust/low-thrust tradeoff is the concept of *multiarcs* between the nodes of a network. This was first proposed by Ishimatsu et al. [26] and subsequently used in state-of-the-art space logistics methods [25, 27-30] to model several space transportation options between the nodes of the space network. In this paper, multiarcs are used to not only represent the transportation of commodities within different types of space vehicles but also to model the different servicer trajectory options available to the optimizer. For example, a servicer may have the option to fly on a 10-day or a 20-day rendezvous maneuver depending on the time it has left to reach and serve a customer satellite. Similarly, a multimodal servicer can fly along a high- or low-thrust trajectory depending on the propulsion system it decides to use. The optimizer's choice of the transportation arcs in the final solution is ultimately driven by the tradeoff between servicing responsiveness and operating costs.

To summarize, this paper generalizes the state-of-the-art OOS logistics method by enabling the automated tradeoff between high-thrust, low-thrust, and multimodal servicers. This is achieved by (1) combining the concepts of dynamic generalized multi-commodity network and piecewise linear approximation for the non-linear low-thrust trajectory model; and (2) leveraging the concept of *multiarcs* to effectively give the optimizer discrete user-defined options related to the propulsion systems and trajectories of the servicers.

Although the developed method can be applied to various OOS concepts, we consider the following OOS concept in this paper as an example. A servicing company deploys high-thrust, low-thrust, and/or multimodal servicers, as well as orbital depots that store the commodities needed to support the operations of the servicers (*e.g.*, propellant, spares). The servicers are equipped with robotic tools designed to perform the operations required for a specific set of services. Launch vehicles may resupply the depots and/or servicers with new commodities if needed. Whenever a GEO satellite requires a service, the OOS operator first decides whether to provide it, and, if so, dispatches the adequate servicer. After performing its task, the servicer then either comes back to its orbital staging location for storage, or flies to another customer satellite if need be. The notional services modeled in this paper are inspection, refueling, station keeping, satellite repositioning, repair, mechanism deployment, and retirement. They fall in either of two categories: random (*e.g.*, unplanned such as component failure), or deterministic (*e.g.*, pre-planned such as refueling).

The rest of this paper is organized as follows. Section II presents the network, the modeling of the servicers' flights, and the generalized MILP model capturing the tradeoff between the high- and low-thrust propulsion alternatives. Section III gives the low-thrust trajectory model currently implemented in the framework. The high-thrust trajectory model used in this paper can be found in Ref. 22. Section IV demonstrates the value of the proposed OOS logistics framework through the applications of short-term operational scheduling and long-term strategic planning of sustainable OOS infrastructures. Finally, section V concludes the paper.

## II. On-Orbit Servicing Logistics Formulation

In this section, the state-of-the-art OOS logistics formulation [22] is generalized by enabling the modeling of high-thrust, low-thrust, and multimodal servicers. Subsection II.A presents the static and dynamic networks used to represent the operations of the OOS infrastructures and the flights of the servicers. Subsection II.B describes how the servicers' trajectory models are computed by user-defined *plugins* interfacing with the proposed OOS logistics framework. Subsection II.C gives the new MILP model which is able to capture the tradeoff between high- and low-thrust propulsion technologies.

### A. Network

In order to properly define the MILP model, the network is developed in two steps. First, the static network is defined, and important modeling elements are introduced. Then, the static network is expanded over time, thus resulting in a dynamic network upon which the MILP model is built.

#### i. Static Network

Figure 1 gives an overview of the static network. This paper models three types of nodes in the network. The *Earth Nodes* are the source nodes of the network from which the commodities are launched into space. The *OOS Parking Nodes* are the orbital slots where depots may be permanently staged, and servicers are parked when idle. The *Customer Nodes* are the orbital slots where the GEO satellites are located and where the service needs are generated. Note that the Orbital Nodes (*i.e.*, OOS Parking Nodes and Customer Nodes) are all located on the same circular orbit which is sufficient to model GEO servicing operations. Later iterations of this framework will relax this assumption to model the operations of OOS infrastructures targeted at other servicing markets than the GEO market.

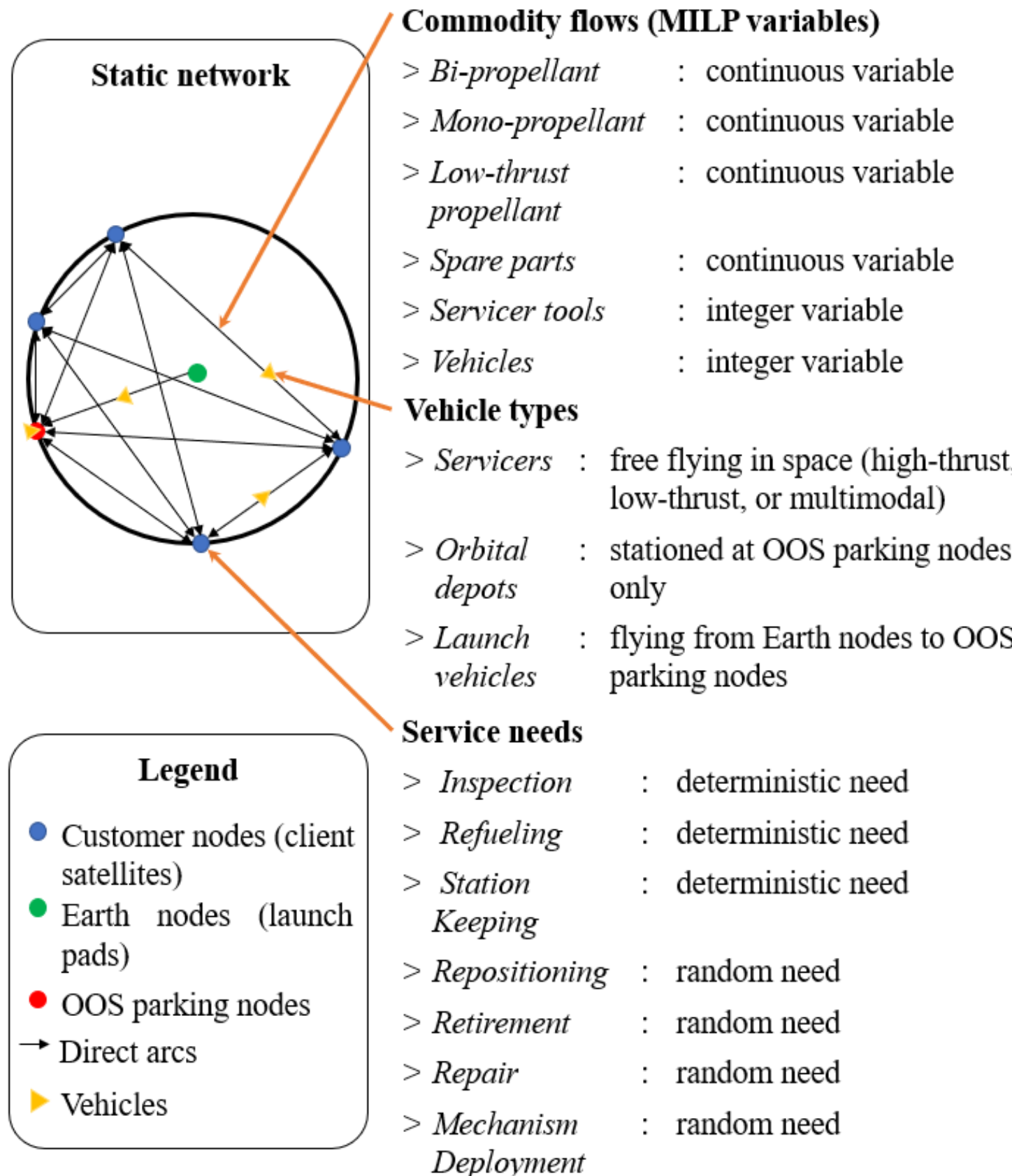

**Fig. 1 High-level overview of the static OOS logistics network; modified from Ref. [22].**

Once the nodes of the network are defined, the transportation arcs between any two nodes are drawn. Multiarcs may be drawn for the same pair of nodes to represent the transportation of commodities with different types of vehicles and over different trajectories. Three types of vehicles are modeled in this paper: (1) the launch vehicles transport commodities from the Earth Nodes to the OOS Parking Nodes; (2) the depots are fixed at the OOS Parking Nodes and store the commodities needed to support the long-term operations of the servicers; (3) the servicers travel between the Orbital Nodes and provide services at the Customer Nodes. Note that the users of the framework can define and assess as many depot and servicer designs as desired, in which case each design is assigned an individual transportation arc for each pair of nodes in the network.

The service needs modeled in the framework occur on a deterministic or random basis. The deterministic service needs comprise Inspection, Refueling, and Station keeping. The random service needs comprise Repositioning, Retirement, Repair, and Mechanism Deployment. The definitions of each of these services are adopted from Ref. [22]. The main parameters used to generate the service needs can also be found in Appendix in Table A1.

The commodities considered in this paper are the propellant needed by the servicers and customer satellites, the spare parts required to perform the repairs, the vehicles, and the tools of the servicers. Bi-propellant is used by the high-thrust and multimodal servicers while the low-thrust propellant is consumed by the low-thrust and multimodal servicers. The mono-propellant is used exclusively for the station keeping of the customer satellites. The spare parts are undifferentiated and thus constitute a continuous commodity, but this can easily be modified to simulate more realistic scenarios. Finally, four servicer tools are modeled: refueling apparatus (T1); Observation sensors (T2); dexterous robotic arm (T3); and capture mechanism (T4). Each tool is used to provide a subset of the services defined previously. The mapping from each tool to a specific subset of the service needs is reminded in Table A2 found in Appendix.

### ii. Dynamic Network

The static network given in Fig. 1 is expanded at predefined time steps by following a periodic topology as illustrated in Fig. 2 with a simple 3-node network. The service time interval of the network, denoted $T$, is chosen to be at most as large as the greatest common divisor of the durations of the considered services. The network is replicated at every $T$ and at two additional time steps within each $T$. This essentially creates 3 intervals per service time interval $T$. The first two are used to enable short spaceflights such as that of the high-thrust servicers and launch vehicles. The third one is used for the provision of services. Figure 2 illustrates the route of a high-thrust servicer in bold yellow arrows. It is sent to orbit by a launch vehicle between $t = 0$ and $t = 2$ before performing on its own an orbital maneuver to rendezvous with a customer satellite between $t = 2$ and $t = 4$. Once the servicer's task is completed, the servicer goes back to its parking location on a 4-day trajectory between $t = 10$ and $t = 14$.

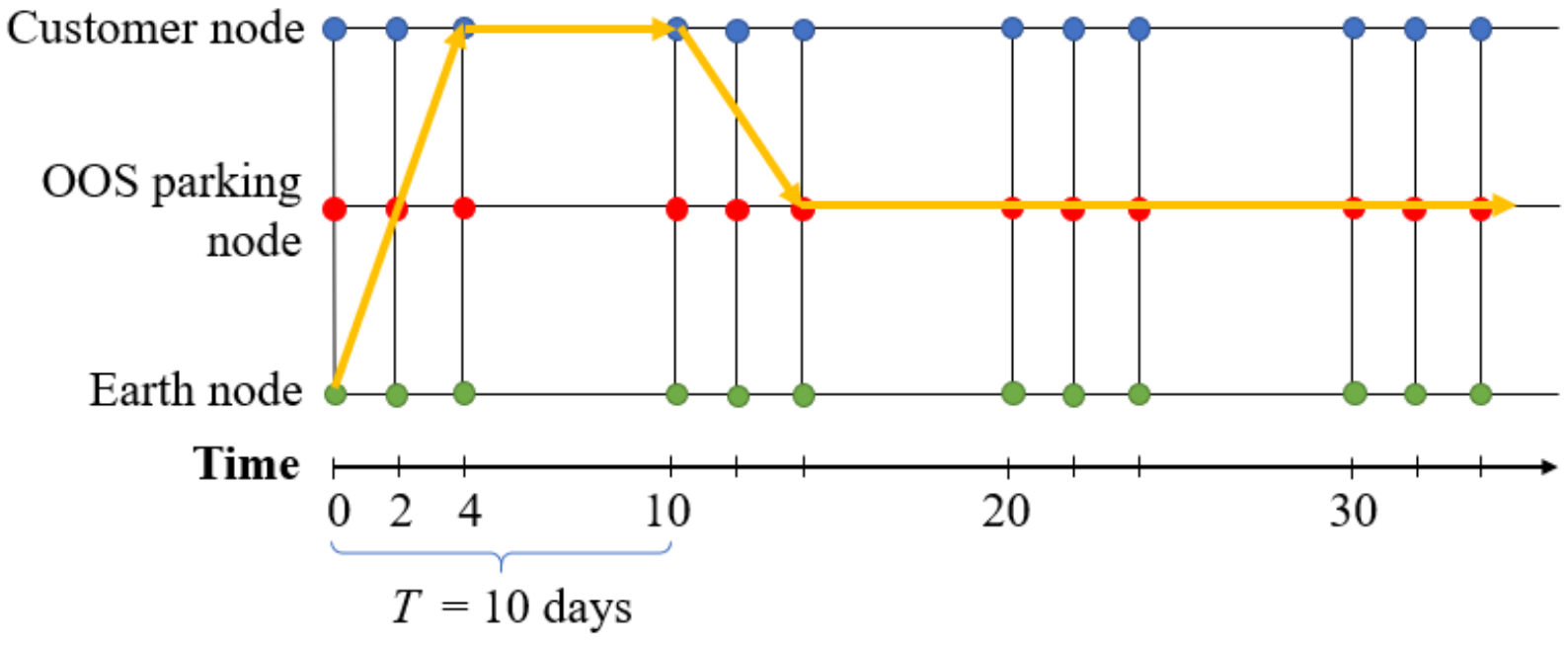

**Fig. 2 Notional dynamic network and high-thrust servicer path; modified from Ref. [22].**

In order to model the long flight durations characteristic of low-thrust servicers and trajectories, the transportation arcs are allowed to span more than the first two intervals defined in each $T$. This is depicted in Fig. 3 with the notional path of a low-thrust servicer. The servicer is sent to orbit on a launch vehicle between $t = 0$ and $t = 2$ before performing a 12-day rendezvous maneuver to reach a customer satellite. Once its task is complete, it comes back to its parking location on a 14-day trajectory between $t = 20$ and $t = 34$. Note that the transportation arcs introduced in Fig. 2 and 3 could in fact be multiarcs representing different strategies for the rendezvous maneuvers performed by the servicers. The reader is referred to section III for the definition of the low-thrust rendezvous maneuvers currently implemented in this framework. The high-thrust trajectories are two-impulse phasing maneuvers and can be found in Ref. 22.

Note that in this paper the servicers are chosen to return to the parking node as soon as they become idle. This behavior is so modeled to prevent hazardous servicer traffic around customer satellites and is captured by Eq. 23-24 in subsection II.C.vii. If needed, this can easily be modified in the proposed framework to allow the servicers to stay near a customer satellite while waiting for a subsequent service need anywhere in the network.

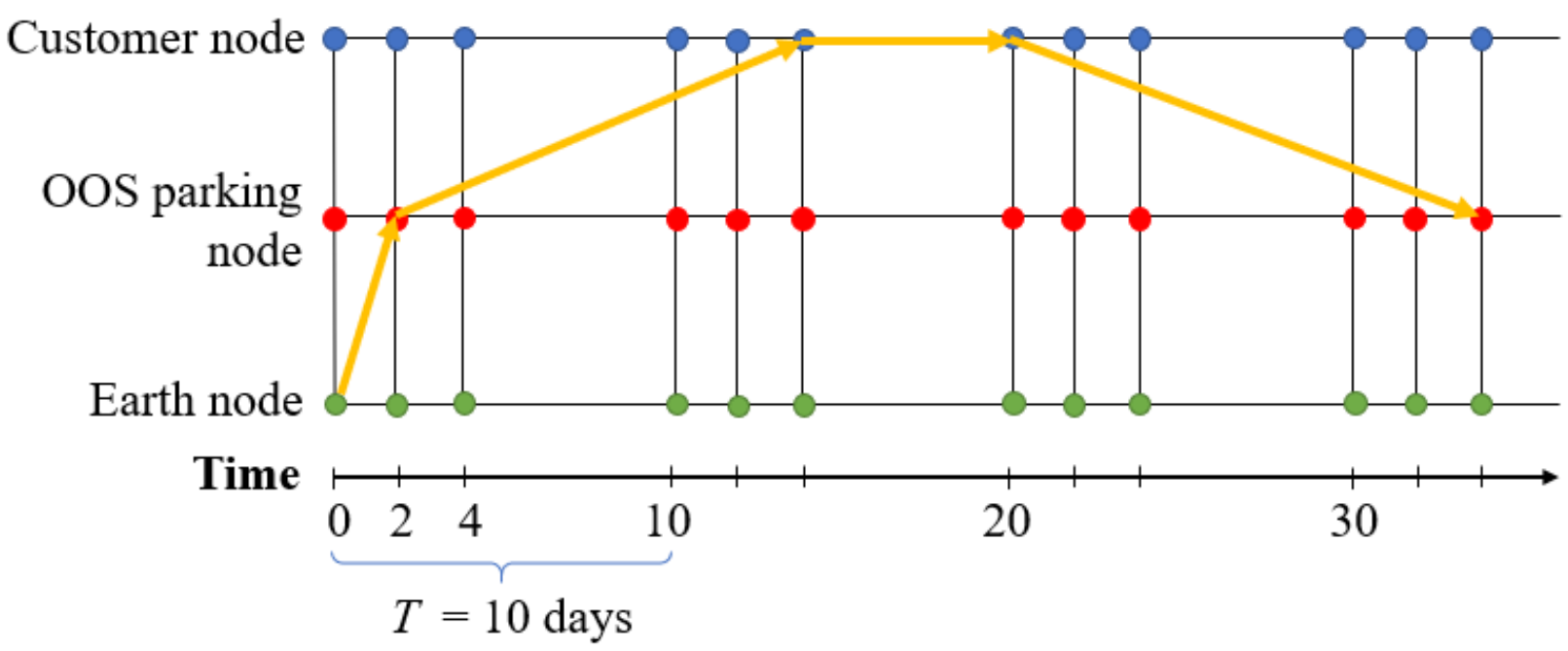

**Fig. 3 Notional path of a low-thrust servicer.**

### B. Servicer Flight Options and Trajectory Plugins

In order to give several servicer flight options to the optimizer, the framework users can define as many flight durations as desired. This allows the optimizer to automatically trade the propellant consumption of the servicers (*i.e.*, cost) and their time of flight (*i.e.*, responsiveness). This is illustrated in Fig. 3 where the servicer first flies from the OOS Parking Node to the Customer Node with a 12-day rendezvous maneuver, and then back to the OOS Parking Node on a 14-day trajectory.

Similarly, the optimizer can be given multiple trajectory options for the flights of the servicers. For example, a high-thrust servicer could use a 2-impulse maneuver or a 3-impulse maneuver to rendezvous with a customer satellite. The optimizer would then choose which trajectory is better suited to maximize the profits generated by the overall OOS infrastructure. The trajectory options given to the optimizer are modeled in user-defined trajectory plugins interfacing with the OOS logistics framework through standard inputs and outputs described in Table 1. Each plugin first computes a mathematical model of the propellant consumption of a servicer as a function of its initial mass. This model is computed over a finite range of the initial mass of the servicer due to the varying masses of propellant and spares in the servicers' tanks and cargo bays. If the model is non-linear – often true of low-thrust models – it cannot be used as is in a MILP. In this case, the plugin automatically converts the non-linear model into a piecewise linear model and is incorporated into the MILP formulation as introduced in the next subsection. In addition, the plugin computes the initial mass of the servicer beyond which the trajectory becomes infeasible (referred to as *servicer mass upper bound* in this paper). The servicer mass upper bound is critical to ensure that overloaded servicers cannot fly over infeasible trajectories (*cf* Eq. 26).

A typical trajectory plugin workflow is summarized in Fig. 4. The high-thrust and low-thrust trajectory models are each integrated into a plugin interfacing with the OOS logistics framework. As many trajectory plugins can be defined and interfaced with the OOS logistics framework as long as they accept the standard inputs and yield the standard outputs given in Table 1.

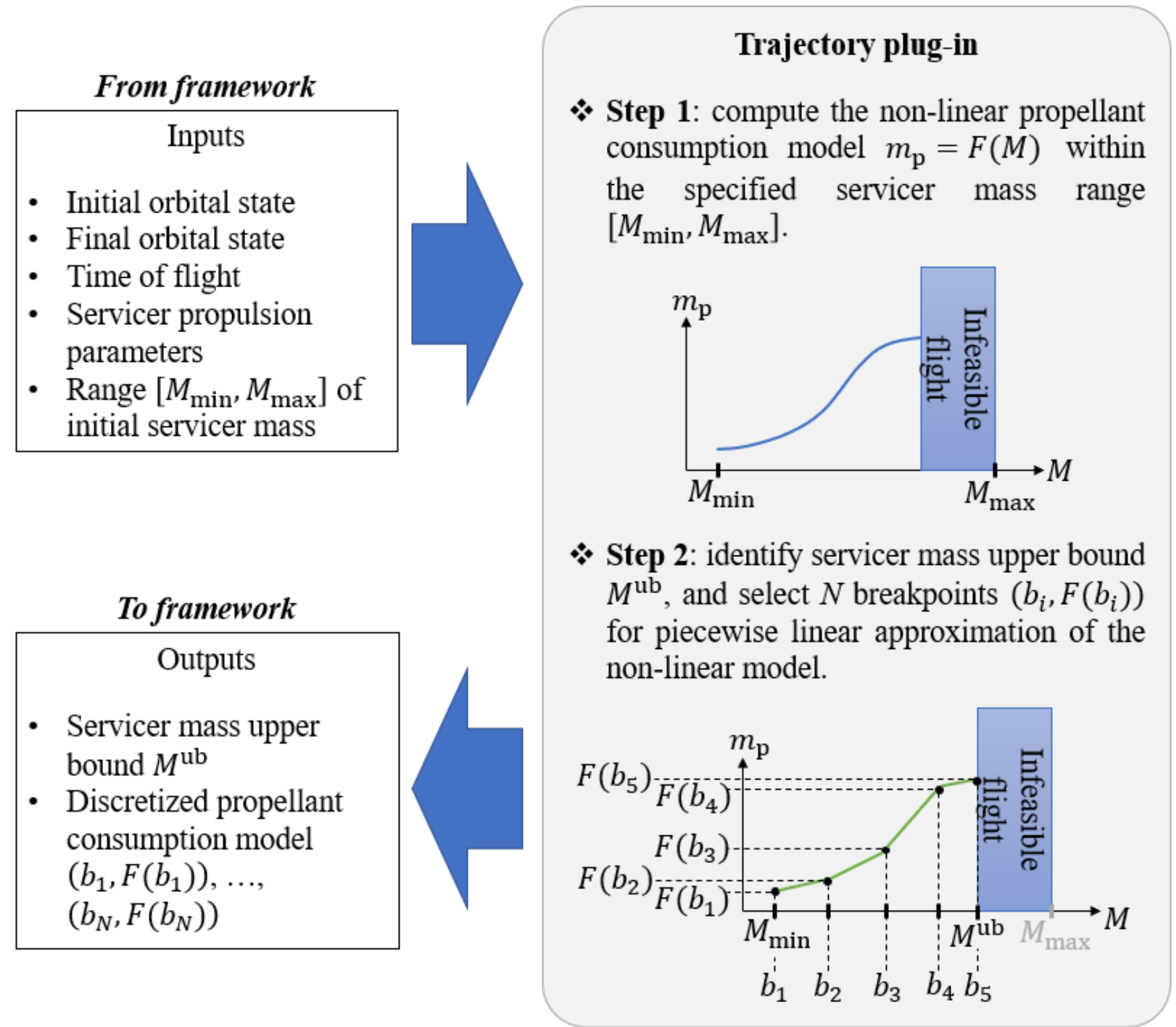

**Fig. 4 Typical trajectory plugin workflow and interfacing with OOS logistics framework.**

**Table 1. Definitions of the trajectory plugins' standardized input and output parameters.**

| Parameter | Definition |
| --- | --- |

| | | |
|---|---|---|
| **Inputs** | *Initial orbital state* | Initial position of the servicer in the network at the beginning of the rendezvous maneuver. |
| | *Final orbital state* | Final position of the servicer in the network at the end of the rendezvous maneuver. |
| | *Time of flight* | Duration of the rendezvous maneuver. |
| | *Servicer propulsion parameters* | Parameters needed to compute the propellant consumption model, such as the specific impulse and the thrust of the engine. |
| | *Range of initial servicer mass* ($[M_{min}, M_{max}]$ in Fig. 4) | Bounds of the mass of the servicer within which the propellant consumption model is computed. The trajectory plugin computes the propellant consumption model strictly within that mass range. |
| **Outputs** | *Servicer mass upper bound* ($M^{ub}$ in Fig. 4) | Initial mass of the servicer beyond which the trajectory becomes infeasible. This may happen in the case of a low-thrust servicer trajectory if the specified time of flight is too short. This parameter is computed by the trajectory plugin. |
| | *Discretized propellant consumption model* | Breakpoints selected from the computed non-linear propellant consumption model for piecewise linear approximation and integration within a MILP. A breakpoint corresponds to an abscissa-ordinate pair, with the abscissa being the initial mass of the servicer and the ordinate being the associated propellant mass consumed. |

### C. Mathematical Formulation

This subsection gives the complete MILP model used to optimize the operations of sustainable OOS infrastructures. This model is built upon the dynamic network presented in section II.A.ii. The state-of-the-art OOS logistics MILP formulation is generalized in three ways:

(1) New index sets and variables are introduced in II.C.i to model the multiarcs associated with the various propulsion system and trajectory options of the servicers;

(2) Additional constraints are defined to properly model the low-thrust rendezvous maneuvers of low-thrust and multimodal servicers (*cf* Eq. 17-19, 26);

(3) Constraints are added to model realistic servicer operations, whether they be related to the provision of services (*cf* Eq. 24) or to the flights of the servicers (*cf* Eq. 27).

After the model is formulated, it is solved either once to optimize the short-term scheduling of the OOS operations, or many times in a row, using the RH procedure, to assess the value of a candidate OOS architecture over the long term. RH decision making is a business practice commonly used in dynamic stochastic settings by which the most immediate decisions are made based on a forecast of near- to medium-term relevant information (*e.g.,* demand for a product or service).

### i. Index Sets and Variables

Consider a time-expanded network graph $\mathcal{G}$ comprised of a set of nodes $\mathcal{N}$ and a set of directed arcs drawn between any two nodes. The set of nodes $\mathcal{N}$ is the union of 3 subsets: the set $\mathcal{E}$ of Earth Nodes; the set $\mathcal{P}$ of OOS Parking Nodes; and the set $\mathcal{C}$ of Customer Nodes. For convenience in writing the objective function and constraints, a set of *Orbital Nodes* is defined as $\mathcal{N}_{\text{orb}} = \mathcal{P} \cup \mathcal{C}$. The nodes of the network are indexed by the letters $i$ and $j$.

There are two types of directed arcs in the network: holdover and transportation arcs. A holdover arc connects a node with itself at two successive time steps. A transportation arc connects two different nodes at two different, not necessarily successive, time steps. Examples of holdover and transportation arcs are illustrated in Fig. 2 and 3 as thin blue lines and bold yellow arrows respectively. Holdover arcs are characterized by three indices: the type of vehicle $v$ flying along the arc, a node $i$ which is both the start and end node of the arc, and the start time $t$ of the arc. Transportation arcs are characterized by six indices: the type of vehicle $v$ flying along the arc, the start node $i$, the end node $j \neq i$, the time of flight $q$ of the rendezvous maneuver, the type of trajectory $r$, and the start time $t$ of the arc. Note that the indices $q$ and $r$ generalize the state-of-the-art OOS logistics MILP formulation by enabling multiarcs that represent various propulsion and trajectory options.

The vehicle types $v$ are defined in index set $\mathcal{V}$. $\mathcal{V}_{\text{serv}}$ and $\mathcal{V}_{\text{dep}}$ are subsets of $\mathcal{V}$ which encompass the different servicer designs and depot designs, respectively. A set $\mathcal{Q}_v$ of time of flight options (index: $q$) is defined per vehicle type $v$. This allows the users of the framework to define a different set of flight duration values per servicer design. Similarly, the set $\mathcal{R}_{vq}$ (index: $r$) stores the trajectory options per vehicle type $v$ and time of flight $q$. The index $r$ thus refers to the trajectory models defined in the user-defined plugins. Finally, the set of time steps $t$ is denoted by $\mathcal{T}$.

The commodities are represented with the index $k$ taking its values in the set $\mathcal{K}$. Subsets of $\mathcal{K}$ are the set $\mathcal{K}_{\text{cont}}$ of continuous commodities (*e.g.*, propellant), and the set $\mathcal{K}_{\text{int}}$ of integer commodities. A subset of integer commodities is the set of servicer tools denoted with $\mathcal{K}_{\text{tools}}$. The flow of commodities along the arc is captured by four sets of variables. The variables $X^{\pm}_{vitk}$ represent the node outflow (superscript +) and node inflow (superscript -) of commodity $k$ along holdover arcs. They are continuous or integer nonnegative variables based on the nature of commodity $k$. Similarly, the binary variables $Y^{\pm}_{vit}$ represent the node outflow and inflow of the vehicle commodity along holdover arcs. The variables $U^{\pm}_{vijqrtk}$ represent the node outflow and inflow of commodity $k$ along transportation arcs. They are continuous or integer nonnegative variables based on the nature of commodity $k$. Finally, the binary variables $W_{vijqrt}^{\pm}$ represent the node outflow and inflow of the vehicle commodity along transportation arcs. For conciseness in certain constraints, the non-vehicle commodity variables are represented under their column vector form as $\boldsymbol{X}_{vit}^{\pm}$ and $\boldsymbol{U}_{vijqrt}^{\pm}$ along holdover arcs and transportation arcs, respectively.

Lastly, the management of services is captured through two additional sets of binary variables. The first one is the set of *service assignment* variables $H_{vst}$ which specify at what time step $t$ the servicer $v \in \mathcal{V}_{\text{serv}}$ must start addressing service need $s$. The service need is addressed if $H_{vst} = 1$. These variables are defined over a subset of $\mathcal{T}$, denoted $\mathcal{W}_s$ (indexed by service need $s$), that encompasses the optional dates to start addressing service need $s$. $\mathcal{W}_s$ is referred to as the *service window* within which the service must start if the optimizer decides to provide it.

The second set of service management variables is the set of *servicer dispatch* variables $B_{vst}$ which are needed to ensure the presence of the servicers at the customer satellites where they are assigned to by the optimizer. The idea behind these variables is that setting $B_{vst}$ to 1 should enforce the presence of servicer $v$ at time step $t$ at whatever Customer Node that triggers the need $s$. Unlike the service assignment variables, these variables are defined for the entire set of time steps $\mathcal{T}$.

Finally, let's define the set $\mathcal{S}_i$ of service needs occurring at customer satellite $i \in \mathcal{C}$, and the set of all service needs $\mathcal{S} = \bigcup_{i \in \mathcal{C}} \mathcal{S}_i$. These sets are indexed by the letter $s$.

**ii. OOS Infrastructure Operation Assumptions**

Before moving on to the mathematical formulation, we first introduce here the assumptions based on which the objective function and constraints are defined:

(1) The launch vehicles are allowed to fly only between the Earth Nodes and the OOS Parking Nodes. This also means that they cannot be staged along the holdover arcs defined at the OOS Parking Nodes. This is to keep the optimizer from leveraging the launch vehicles as makeshift depots;

(2) The orbital depots remain staged at the OOS Parking Nodes;

(3) One orbital depot is deployed at every OOS Parking Node defined in the network. This allows us to model sustainable OOS infrastructures that support the long-term operations of servicers (*e.g.*, through regular servicer refueling);

(4) Each orbital depot is assumed to be a small-scale robotic space station technically capable of operating the transfer of non-vehicle commodities across the different spacecraft (*e.g.*, from the fairing of a launch vehicle to a servicer's cargo bay);

(5) The servicers are not capable of exchanging commodities with each other. They may exchange non-vehicle commodities at the OOS Parking Nodes through the orbital depots;

(6) There cannot be more than one servicer at a time at a Customer Node. This prevents undesirable traffic of servicers near customer satellites;

(7) A servicer can provide only one service at a time to a single customer satellite, even if it is equipped with several tools;

(8) A servicer cannot depart the OOS Parking Node at which it is staged unless it has been assigned a service need;

(9) A servicer must go back to one of the available OOS Parking Nodes after providing a service unless it has been assigned a sequence of back-to-back services;

(10) The operating costs of the servicers are the same no matter whether they are idle at a node, flying between two different nodes, or providing a service to a customer satellite;

(11) No penalty fee is incurred to the OOS infrastructure if the optimizer decides not to address a service need. Penalty fees, however, are incurred if the optimizer decides to provide a service and does so with delays. Here, the delay is defined as the difference between the time step at which a service need is triggered and the time step at which the service starts.

**iii. Objective Function**

This paper assumes that the goal of OOS businesses is to maximize their profits, defined as the revenues generated by the provision of services minus the costs to manufacture, deploy and operate the OOS infrastructures. The revenues in the MILP model are captured by Eq. 1,

$$J_{\text{revenues}} = \sum_{v \in \mathcal{V}_{\text{serv}}} \sum_{i \in \mathcal{C}} \sum_{s \in \mathcal{S}_i} \sum_{t \in \mathcal{W}_s} \{ r_s H_{vst} \} \tag{1}$$

where $r_s$ is the revenue generated if service need $s$ is taken care of by the OOS infrastructure.

Equation 2 gives the expression for the launch cost,

$$J_{\text{launch}} = c^{\text{launch}} \sum_{v \in \mathcal{V}} \sum_{t \in \mathcal{T}} \sum_{i \in \mathcal{E}} \sum_{j \in \mathcal{P}} \left\{ \sum_{k \in \mathcal{K}_{\text{cont}}} U^{+}_{vijq_0r_0tk} + \sum_{k \in \mathcal{K}_{\text{int}}} m_k U^{+}_{vijq_0r_0tk} + m_v W^{+}_{vijq_0r_0t} \right\} \tag{2}$$

where $m_k$ and $m_v$ are the mass per unit of integer commodity $k$ and vehicle $v$ respectively, $c^{\text{launch}}$ is the launch cost per unit mass, and $q_0$ and $r_0$ are special indices used to characterize the launches from the Earth Nodes to the OOS Parking Nodes.

Equation 3 expresses the cost associated with developing and manufacturing the servicing elements, and purchasing the consumables to support the OOS infrastructure,

$$J_{\text{pdm}} = \sum_{v\in\mathcal{V}} \sum_{t\in\mathcal{T}} \sum_{i\in\mathcal{E}} \sum_{j\in\mathcal{P}} \left\{ \sum_{k\in\mathcal{K}} c_k^{\text{pdm}} U_{vijq_0r_0tk}^{+} + c_v^{\text{pdm}} W_{vijq_0r_0t}^{+} \right\} \tag{3}$$

where $c_k^{\text{pdm}}$ is the cost per unit of integer commodity and per unit mass of continuous commodity, $c_v^{\text{pdm}}$ is the cost per unit of vehicle commodity, and the acronym "pdm" stands for *purchase/development/manufacturing costs*.

Equation 4 captures the fees incurred to the OOS infrastructure if a service is provided with delays,

$$J_{\text{delay}} = \sum_{i\in\mathcal{C}} \sum_{s\in\mathcal{S}_i} \sum_{\tau\in\mathcal{W}_s} \left\{ c_s^{\text{delay}} (\tau - \tau_s) \sum_{v\in\mathcal{V}_{\text{serv}}} H_{vs\tau} \right\} \tag{4}$$

where $c_s^{\text{delay}}$ is the fee per unit time of delay, and $\tau_s$ is the time step representing the beginning of the service window associated with service need $s$. Note that $\tau_s$ is the smallest of the time steps defined in $\mathcal{W}_s$, so the delay $\tau - \tau_s$ with $\tau \in \mathcal{W}_s$ is nonnegative.

Equation 5 models the operating cost of the orbital depots,

$$J_{\text{dep}} = \sum_{v\in\mathcal{V}_{\text{dep}}} \sum_{i\in\mathcal{P}} \sum_{t\in\mathcal{T}} \left\{ c_v^{\text{dep}} \Delta_t Y_{vit}^{+} \right\} \tag{5}$$

where $c_v^{\text{dep}}$ is the operating cost of the depots per unit time, and $\Delta_t$ is the length of the holdover arc starting at time step $t$.

Finally, Eq. 6 models the operating cost of the servicers,

$$J_{\text{serv}} = \sum_{v\in\mathcal{V}_{\text{serv}}} \sum_{i\in\mathcal{N}_{\text{orb}}} \sum_{t\in\mathcal{T}} c_v^{\text{serv}} \left\{ \Delta_t Y_{vit}^{+} + \sum_{\substack{j\in\mathcal{N}_{\text{orb}} \\ j\neq i}} \sum_{q\in\mathcal{Q}_v} \sum_{r\in\mathcal{R}_{vq}} q W_{vijqrt}^{+} \right\} \tag{6}$$

where $c_v^{\text{serv}}$ is the operating cost of the servicers per unit time.

### iv. Mass Balance Constraints

The mass balance constraints ensure the conservation of commodities at the nodes of the network. Equation 7 first gives the mass balance of non-vehicle commodities at the Customer Nodes. This constraint is written per servicer to model the fact that servicers cannot exchange commodities with each other. Note also that the mass balance constraint at the Customer Nodes is not specialized for launch vehicles and orbital depots as these spacecraft are not allowed to visit the customer satellites.

$\forall v \in \mathcal{V}_{\text{serv}}, \forall i \in \mathcal{C}, \forall t \in \mathcal{T}$:

$$\boldsymbol{X}_{vit}^{+} - \boldsymbol{X}_{vi(t-\Delta_t')}^{-} + \sum_{\substack{j \in \mathcal{N} \\ q \in \mathcal{Q}_v \\ r \in \mathcal{R}_{vq} \\ i \neq j}} \left\{ \boldsymbol{U}_{vijqrt}^{+} - \boldsymbol{U}_{vjiqr(t-q)}^{-} \right\} = \begin{cases} \sum_{s \in \mathcal{S}_i} \boldsymbol{d}_{vs} H_{vst} & \text{if } t \in \mathcal{L}_i \\ 0 & \text{otherwise} \end{cases} \tag{7}$$

where $\Delta_t'$ is the length of the holdover arc directly preceding time step $t$, $\boldsymbol{d}_{vs}$ is the column vector of nonpositive demand for the non-vehicle commodities (*e.g.*, propellant for refueling), and $\mathcal{L}_i$ is the set of time steps which fall within at least one of the service windows $\mathcal{W}_s$ defined for the service needs $s$ occurring at customer satellite $i$. The set $\mathcal{L}_i$ is defined in Eq. 8:

$$\mathcal{L}_i = \bigcup_{s \in \mathcal{S}_i} \mathcal{W}_s \tag{8}$$

Equation 9 gives the mass balance of the non-vehicle commodities at the OOS Parking Nodes. Note that unlike Eq. 7, this equation is written by summing the variables over the vehicle index $v \in \mathcal{V}$ to allow the spacecraft to exchange commodities with each other through the orbital depots,

$\forall i \in \mathcal{P}, \forall t \in \mathcal{T}$:

$$\sum_{v \in \mathcal{V}} \left\{ \boldsymbol{X}_{vit}^{+} - \boldsymbol{X}_{vi(t-\Delta_t')}^{-} \right\} + \sum_{\substack{v \in \mathcal{V} \\ j \in \mathcal{N} \\ q \in \mathcal{Q}_v \\ r \in \mathcal{R}_{vq} \\ i \neq j}} \left\{ \boldsymbol{U}_{vijqrt}^{+} - \boldsymbol{U}_{vjiqr(t-q)}^{-} \right\} = 0 \tag{9}$$

Equation 10 models the supply of commodities from the Earth Nodes, with $\boldsymbol{\sigma}_{it}$ being the column vector that represents the nonnegative supply of non-vehicle commodities,

$\forall i \in \mathcal{E}, \forall t \in \mathcal{T}$:

$$\sum_{v \in \mathcal{V}} \left\{ \boldsymbol{X}_{vit}^{+} - \boldsymbol{X}_{vi(t-\Delta_t')}^{-} \right\} + \sum_{\substack{v \in \mathcal{V} \\ j \in \mathcal{N} \\ q \in \mathcal{Q}_v \\ r \in \mathcal{R}_{vq} \\ i \neq j}} \left\{ \boldsymbol{U}_{vijqrt}^{+} - \boldsymbol{U}_{vjiqr(t-q)}^{-} \right\} \leq \boldsymbol{\sigma}_{it} \tag{10}$$

Finally, Eq. 11 and 12 represent the mass balance of the non-vehicle commodity for each type of spacecraft at the Orbital Nodes (*i.e.*, Customer Nodes and OOS Parking Nodes) and Earth Nodes, respectively,

$\forall v \in \mathcal{V}, \forall i \in \mathcal{N}_{\text{orb}}, \forall t \in \mathcal{T}$:

$$Y_{vit}^{+} - Y_{vi(t-\Delta_t')}^{-} + \sum_{\substack{j \in \mathcal{N} \\ q \in \mathcal{Q}_v \\ r \in \mathcal{R}_{vq} \\ i \neq j}} \left\{ W_{vijqrt}^{+} - W_{vjiqr(t-q)}^{-} \right\} = 0 \tag{11}$$

$\forall v \in \mathcal{V}, \forall i \in \mathcal{E}, \forall t \in \mathcal{T}$:

$$Y_{vit}^{+} - Y_{vi(t-\Delta_t')}^{-} + \sum_{\substack{j \in \mathcal{N} \\ q \in \mathcal{Q}_v \\ r \in \mathcal{R}_{vq} \\ i \neq j}} \left\{ W_{vijqrt}^{+} - W_{vjiqr(t-q)}^{-} \right\} \leq \sigma_{ivt}' \tag{12}$$

where $\sigma_{ivt}'$ represents the nonnegative supply of spacecraft at the Earth Nodes.

### v. Concurrency Constraints

The concurrency constraints relate to the payload capacities of the vehicles considered in the network. This is to enforce the transportation of commodities (*e.g.*, propellant, spares) within vehicles. Equation 13 represents the concurrency constraints along holdover arcs, while Eq. 14 is for the concurrency constraints along transportation arcs.

$\forall v \in \mathcal{V}, \forall i \in \mathcal{N}_{\text{orb}}, \forall t \in \mathcal{T}$:

$$\bar{\bar{M}}_{vii} \boldsymbol{X}_{vit}^{+} \leq \boldsymbol{e}_v Y_{vit}^{+} \tag{13}$$

$\forall v \in \mathcal{V}, \forall i \in \mathcal{N}, \forall j \in \mathcal{N}, i \neq j, \forall q \in \mathcal{Q}_v, \forall r \in \mathcal{R}_{vq}, \forall t \in \mathcal{T}$:

$$\bar{\bar{M}}_{vij} \boldsymbol{U}_{vijqrt}^{+} \leq \boldsymbol{e}_v W_{vijqrt}^{+} \tag{14}$$

where $\bar{\bar{M}}_{vij}$ is the concurrency matrix, and $\boldsymbol{e}_v$ represents the vehicle design parameters such as payload and propellant capacities.

### vi. Commodity Transformation Constraints

These constraints model the consumption of the commodities needed to operate the depots and servicers (*e.g.*, propellant consumption). Equations 15 and 16 give the constraints for holdover and transportation arcs respectively when the consumption models are linear.

$\forall v \in \mathcal{V}, \forall i \in \mathcal{N}, \forall t \in \mathcal{T}$:

$$\bar{\bar{Q}}_{vi}' \begin{bmatrix} \boldsymbol{X}_{vit}^{+} \\ Y_{vit}^{+} \end{bmatrix} = \begin{bmatrix} \boldsymbol{X}_{vit}^{-} \\ Y_{vit}^{-} \end{bmatrix} \tag{15}$$

$\forall v \in \mathcal{V}, \forall i \in \mathcal{N}, \forall j \in \mathcal{N}, i \neq j, \forall q \in \mathcal{Q}_v, \forall r \in \mathcal{R}_{vq}, \forall t \in \mathcal{T}$: (16)

$$\bar{\bar{Q}}''_{vijqr}\begin{bmatrix}\boldsymbol{U}^{+}_{vijqrt}\\ W^{+}_{vijqrt}\end{bmatrix}=\begin{bmatrix}\boldsymbol{U}^{-}_{vijqrt}\\ W^{-}_{vijqrt}\end{bmatrix}$$

where $\bar{\bar{Q}}'_{vi}$ and $\bar{\bar{Q}}''_{vijqr}$ are the mass transformation matrices for holdover and transportation arcs, respectively.

When the consumption models are non-linear, as is often the case of low-thrust models, Eq. 15 and 16 cannot properly capture their non-linear nature. Instead, they must be converted into piecewise linear approximations, and as many additional variables as the number of breakpoints must be introduced. More specifically in this paper, the model describing the propellant consumption of a servicer with respect to its total initial mass is approximated and integrated into the MILP. We now define $P^{\pm}_{vijqrt}$ as the inflow/outflow variables of the propellant meant to be consumed by the servicers along the transportation arcs. As depicted in Fig. 4, we also denote by $F(.)$ the considered non-linear propellant consumption model, and by $(b_1, F(b_1))$, …, $(b_N, F(b_N))$ the $N$ breakpoints used to define the piecewise linear approximation. An SOS2* set of $N$ continuous nonnegative variables $\lambda_1$, …, $\lambda_N$ is then introduced, and the approximated value of the propellant consumed over the transportation arc (*i.e.*, $P^{+}_{vijqrt} - P^{-}_{vijqrt}$) can be found as [31]:

$\forall v \in \mathcal{V}_{\text{serv}}, \forall i \in \mathcal{N}, \forall j \in \mathcal{N}, i \neq j, \forall q \in \mathcal{Q}_v, \forall r \in \mathcal{R}_{vq}, \forall t \in \mathcal{T}$:

$$\lambda_1 + \cdots + \lambda_N = 1 \tag{17}$$

$$\lambda_1 b_1 + \cdots + \lambda_N b_N = Z^{+}_{vijqrt} \tag{18}$$

$$\lambda_1 F(b_1) + \cdots + \lambda_N F(b_N) = P^{+}_{vijqrt} - P^{-}_{vijqrt} \tag{19}$$

where $Z^{+}_{vijqrt}$ is the total mass of the servicer along the transportation arc as defined in Eq. 20. Note that this model is a newly developed component to the OOS logistics framework, as it was not needed when the high-thrust linear rocket equation model was the only trajectory option considered in the past work.

$\forall v \in \mathcal{V}_{\text{serv}}, \forall i \in \mathcal{N}, \forall j \in \mathcal{N}, i \neq j, \forall q \in \mathcal{Q}_v, \forall r \in \mathcal{R}_{vq}, \forall t \in \mathcal{T}$:

$$Z^{+}_{vijqrt} = \sum_{k \in \mathcal{K}_{\text{cont}}} U^{+}_{vijqrtk} + \sum_{k \in \mathcal{K}_{\text{int}}} m_k U^{+}_{vijqrtk} + m_v W^{+}_{vijqrt} \tag{20}$$

* An SOS2 (special ordered sets of type 2) constraint enforces that, at most, two of the $\lambda$ can be nonzero; and these two nonzero elements of the set must be consecutive.

#### vii. Service Management Constraints

Five constraints are defined for the optimizer to properly manage the services and the associated operations of the fleet of servicers. Equation 21 models the assignment of the service needs to the servicers. It ensures that each service need $s$ is assigned at most once to a servicer within the service window $\mathcal{W}_s$.

$\forall s \in \mathcal{S}$:

$$\sum_{v \in \mathcal{V}_{\text{serv}}} \sum_{\tau \in \mathcal{W}_s} H_{vs\tau} \leq 1 \tag{21}$$

Equation 22 describes the coupling between the service assignment variables and servicer dispatch variables through the binary parameters $\beta_{s\tau t}$. These parameters are automatically generated by the OOS optimization framework before running the optimization based on the input data related to the service needs. Note that these parameters are indexed by both $t \in \mathcal{T}$ and $\tau \in \mathcal{W}_s$. Essentially, the binary parameter $\beta_{s\tau t}$ captures the duration of service $s$ and at which time steps $t$ the service has to be provided if the OOS operator decides to start the service at the date $\tau$.

$\forall v \in \mathcal{V}_{\text{serv}}, \forall s \in \mathcal{S}, \forall t \in \mathcal{T}$:

$$B_{vst} = \sum_{\tau \in \mathcal{W}_s} H_{vs\tau} \beta_{s\tau t} \tag{22}$$

Equation 23 is designed to keep the optimizer from providing more than one service at a time to a customer satellite. One and only one service need can be addressed at a time at each Customer Node of the network. This constraint along with Eq. 24 also avoids the unnecessary clustering of servicers at a Customer Node for the purpose of increased operational safety. Indeed, GEO satellite operators are likely to require minimal servicer traffic around their spacecraft.

$\forall i \in \mathcal{C}, \forall t \in \mathcal{T}$:

$$\sum_{v \in \mathcal{V}_{\text{serv}}} \sum_{s \in \mathcal{S}_i} B_{vst} \leq 1 \tag{23}$$

$\forall v \in \mathcal{V}_{\text{serv}}$, $\forall i \in \mathcal{C}$, $\forall t \in \mathcal{T}$:

$$Y_{vit}^{+} = \sum_{s \in \mathcal{S}_i} B_{vst} \tag{24}$$

Finally, Eq. 25 ensures the adequate servicer tool is used to provide a service, as specified by the service-tool mapping given in Table A2 in Appendix. The binary parameter $\gamma_{sk}$ is set to 1 if the tool $k$ is required to address service need $s$, and to 0 otherwise.

$\forall v \in \mathcal{V}_{\text{serv}}, \forall i \in \mathcal{C}, \forall t \in \mathcal{T}, \forall k \in \mathcal{K}_{\text{tools}}$:

$$X^{+}_{vitk} \geq \sum_{s \in \mathcal{S}_i} \gamma_{sk} B_{vst} \tag{25}$$

### viii. Servicer Spaceflights

Two important constraints are presented here to accurately model the flights of the servicers in a MILP. The first constraint, captured by Eq. 26, ensures that a servicer cannot fly along a transportation arc between any two nodes of the network if its total mass (*i.e.*, combined servicer's dry mass and payload mass) exceeds the servicer mass upper bound calculated by the corresponding trajectory plugin (*cf* subsection II.B).

$\forall v \in \mathcal{V}_{\text{serv}}, \forall i \in \mathcal{N}, \forall j \in \mathcal{N}, i \neq j, \forall q \in \mathcal{Q}_v, \forall r \in \mathcal{R}_{vq}, \forall t \in \mathcal{T}$:

$$Z^{+}_{vijqrt} \leq M^{\text{ub}}_{vijqr} W^{+}_{vijqrt} \tag{26}$$

where $M^{\text{ub}}_{vijqr}$ is the servicer mass upper bound, and $Z^{+}_{vijqrt}$ is the total mass of the servicer and its payloads as defined by Eq. 20.

The second constraint, captured by Eq. 27, makes sure that a servicer departs its parking location only when it is assigned a service need by the optimizer. This is essential for the servicer to go straight from its parking location to its target location without undesirably stopping by other network nodes along the way.

$\forall v \in \mathcal{V}_{\text{serv}}, \forall i \in \mathcal{C}, \forall t \in \mathcal{T}$:

$$\sum_{\substack{j \in \mathcal{N} \\ q \in \mathcal{Q}_v \\ r \in \mathcal{R}_{vq} \\ j \neq i}} W^{+}_{vjiqr(t-q)} = \begin{cases} \sum_{s \in \mathcal{S}_i} H_{vst} & \text{if } t \in \mathcal{L}_i \\ 0 & \text{otherwise} \end{cases} \tag{27}$$

## III. Trajectory Models

The formulation proposed in Section II enables the trajectory models to be used for the OOS mission design optimization. Although the proposed formulation can integrate various trajectory plugins, this section describes an example low-thrust model integrated into a trajectory plugin of the OOS logistics framework. The high-thrust model used in this paper is based on two-impulse phasing maneuvers and can be found in Ref. 22. The procedure needed to derive the propellant consumption as a function of the servicer's initial mass is given. Note that through the standard interface between the OOS logistics framework and the trajectory plugins, the users of the framework can define as many new trajectory models as desired.

## A. Low-Thrust Phasing Maneuver

The analytical low-thrust trajectory model presented in this subsection is based on Ref. [32] (*cf* pp. 5-7, *Re-positioning in Orbits: Walking*). This trajectory is a phasing maneuver for a low-thrust servicer to change its angular position along a circular orbit (*e.g.* GEO orbit). It consists of three phases: (1) an initial thrust phase which brings the servicer from its initial circular orbit to an intermediate circular orbit; (2) a coasting phase along the intermediate circular orbit; and (3) a final thrust phase that brings the servicer back to its initial circular orbit but at a different angular position. Figure 5 illustrates the trajectory. Note that due to the low thrust of the servicer, the initial and final thrust phases are spiraling about the Earth but are not represented as such in Fig. 5 for clarity. The phasing angle $\delta\theta$ defined in Fig. 5 is the angle between the position the servicer would have if it had remained coasting along its initial circular orbit and its actual position with the phasing maneuver. The phasing angle is represented at three different dates $t_1, t_2, t_3$ in Fig. 5.

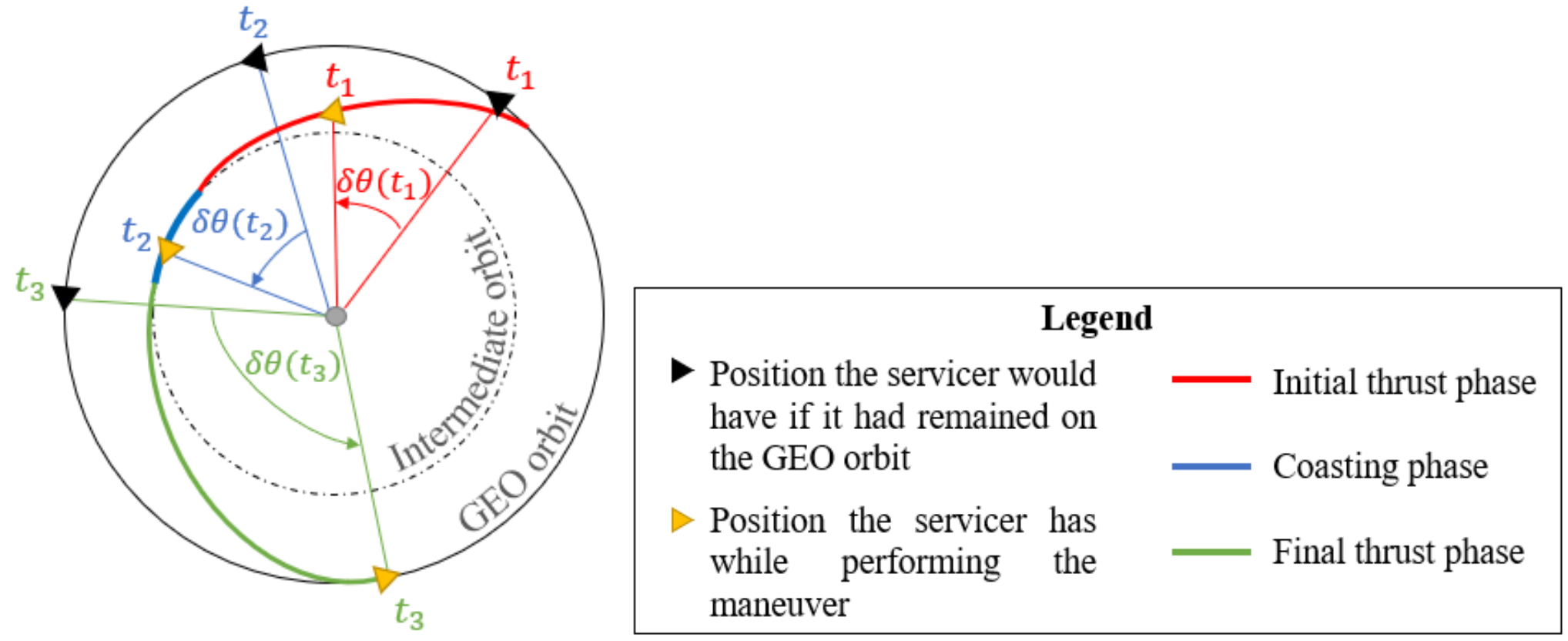

**Fig. 5 Overview of the low-thrust phasing maneuver modeled in this paper.**

A few assumptions are made in Ref. [32] to develop a simple analytical model describing the low-thrust phasing maneuver introduced in Fig. 5:

- The thrust is constant, continuous, and tangential during the thrusting phases;
- The mass of the servicer is assumed constant during the maneuver due to a low propellant consumption. This results in a constant thrust acceleration during the thrusting phases;
- The thrust acceleration is assumed much smaller than the gravitational and centripetal accelerations to approximate the spiraling transfer as a series of circular orbits of varying radii (*spiral approximation*);
- The specific impulse remains constant during the thrusting phases (*i.e.*, the engine runs at a constant power level).

From these assumptions and the power balance, Ref. [32] develops Eq. 28 describing the dynamic evolution of the phasing angle,

$$\frac{d^2(\delta\theta)}{dt^2} = -\frac{3a_\theta}{r_0} \tag{28}$$

where $a_\theta$ is the thrust acceleration due to the engine, and $r_0$ is the radius of the initial circular orbit. Note that the thrust acceleration is backward (*i.e.*, negative) for a positive target phasing angle, and forward (*i.e.*, positive) otherwise. If we note $\Delta\theta$ the target phasing angle defined over (-π,π], then the thrust acceleration is given by Eq. 29 as

$$a_\theta = -\text{sign}(\Delta\theta)\frac{F}{m_0} \tag{29}$$

where $F$ is the thrust of the engine, and $m_0$ is the constant mass of the servicer.

We then integrate Eq. 28 twice for each phase of the maneuver and enforce the continuity in the servicer's phasing angle and rate between the different phases. This leads to Eq. 30, which is a quadratic equation in the duration $\tau$ of the initial and final thrust phases. Note that since the mass of the servicer is assumed constant during the maneuver, the initial and final thrust phases are of the same duration $\tau$.

$$\tau^2 - t_\text{f}\cdot\tau + \frac{r_0 m_0|\Delta\theta|}{3F} = 0 \tag{30}$$

where $t_\text{f}$ is the duration of the entire rendezvous maneuver.

The coefficients of Eq. 30 are known: $\Delta\theta$, $t_\text{f}$, $r_0$ and $F$ are the inputs to the trajectory plugin, and $m_0$ is the initial mass of the servicer, which is varied between some bounds $M_\text{min}$ and $M_\text{max}$ by the trajectory plugin to compute the propellant consumption model associated with the low-thrust phasing maneuver (*cf* II.B and Fig. 4). The *servicer mass upper bound* $M^\text{ub}$ is then found by setting the discriminant of Eq. 30 to 0 and solving for the initial mass. The flight of the servicer is feasible for an initial mass $m_0 \leq M^\text{ub}$ but infeasible otherwise as modeled by Eq. 26.

Assuming that the maneuver is feasible for the given initial mass $m_0$, the trajectory plugin proceeds by computing the propellant consumed during the maneuver assuming a constant propellant mass flow rate $b > 0$,

$$m_\text{p} = 2b\tau \tag{31}$$

From the thrust force $F$ and specific impulse $I_{sp}$ of the low-thrust engine (also an input to the trajectory plugin), the expression for the propellant mass flow rate is:

$$b = \frac{F}{g_0 I_{sp}} \tag{32}$$

The non-linear propellant consumption model is computed and illustrated in Fig. 6 for a constant thrust force $F = 1.16N$, a specific impulse $I_{sp} = 1{,}790s$, a time of flight $\Delta t = 8$ days, a target phasing angle $\Delta\theta = 180°$, and an initial orbit radius $r_0 = 42{,}164$ km corresponding to the GEO orbit. The data for the thrust force and specific impulse are based upon Northrop Grumman's MEV [6,33]. The model is computed within the bounds $M_{\min} = 500$ kg and $M_{\max} = 4000$ kg for the servicer's initial mass. The *servicer mass upper* bound is found to be 3,138 kg as can be seen at the discontinuity in Fig. 6. For a servicer initial mass below 3,138 kg, the maneuver is feasible, but infeasible beyond. In addition to computing this non-linear model, the trajectory plugin would proceed by finding its piecewise linear approximation for proper integration into the MILP.

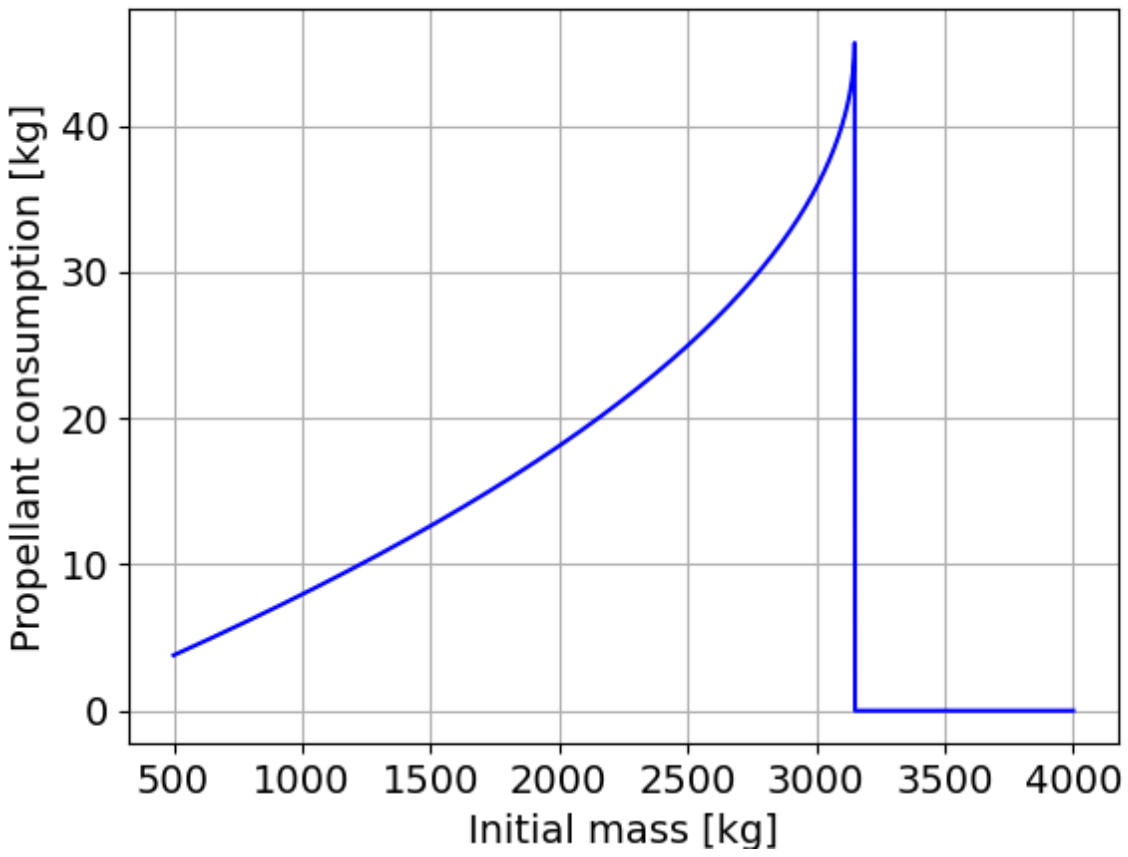

**Fig. 6 Non-linear propellant consumption of the servicer as a function of its initial mass.**

## IV. Case Studies

This section demonstrates two use cases of the OOS logistics optimization framework: a short-term operational scheduling use case; and a long-term strategic planning use case. In both cases, the optimizer now has the capability to automatically trade between the high- and low-thrust propulsion systems and trajectories of the servicers. Subsection IV.A summarizes the assumptions used to run the different scenarios designed to demonstrate the use cases. Subsection IV.B demonstrates the proposed method as a tool for short-term operational scheduling of an existing OOS infrastructure equipped with both high- and low-thrust technologies. Subsection IV.C demonstrates the proposed method as a tool to inform long-term strategic design of OOS architectures, especially regarding high- and low-thrust technologies, based on their performance given various OOS market conditions.

## A. Assumptions and Scenarios

In this subsection, we first present the assumptions and data associated with the fleet of customer satellites. Then, the assumptions and data associated with the OOS infrastructure are given. Finally, the scenarios associated with the two use-cases are presented.

### i. Customer Fleet Assumptions

The data related to the deterministic and random service needs are given in Table 2 and Table 3, respectively; these are the same assumptions used in the case studies for Ref. [22]. Note that the given data are for one customer satellite; by increasing the number of customer satellites in the simulations, the service need rates of the entire fleet of customer satellites increase. The positions of the customer satellites were retrieved from the UCS (Union of Concerned Scientists) database of satellites [34].

**Table 2. Assumptions related to the deterministic service needs (*).**

| Parameters | Inspection | Refueling | Station Keeping |
|---|---|---|---|
| **Revenues [$M]** | 10 [13] | 15 [13] | 20 |
| **Delay penalty fee [$/day]** | 5,000 [13] | 100,000 [13] | 100,000 [13] |
| **Service duration [days]** | 10 [13] | 30 [13] | 180 [19] |
| **Service window [days]** | 30 | 30 | 30 |
| **Frequency of occurrence [days]** | 6,310 [21] | 2,100 [21] | 2,100 [21] |

* References are indicated within brackets. Data without reference is assumed. The frequency of occurrence is derived from the data given in Ref. [21].

**Table 3. Assumptions related to the random service needs (*).**

| Parameters | Repositioning | Retirement | Repair | Mechanism Deployment |
|---|---|---|---|---|
| **Revenues [$M]** | 10 [13] | 10 [13] | 30 | 25 [13] |
| **Delay penalty fee [$/day]** | 100,000 [13] | 0 [13] | 100,000 [13] | 100,000 [13] |
| **Service duration [days]** | 30 [13] | 30 [13] | 30 [13] | 30 [13] |
| **Service window [days]** | 30 | 30 | 30 | 30 |
| **Mean frequency of occurrence [days]** | 2,520 [21] | 2,520 [21] | 9,020 [21] | 21,050 [21] |

* References are indicated within brackets. Data without reference is assumed. The mean frequency of occurrence is derived from the data given in Ref. [21].

### ii. OOS Infrastructure Assumptions

The four notional servicer tools given in Table A2 have an assumed cost of $100,000 and an assumed mass of 100kg. The other commodities considered in the case studies are the spares (assumed price tag: $1,000/kg), bipropellant for the servicers (price tag for Monomethyl Hydrazine: $180/kg), monopropellant (price tag for Hydrazine: $230/kg), and low-thrust propellant (price tag for Xenon gas: $1,115/kg).

An orbital depot is assumed pre-deployed at a GEO orbital slot located at a longitude of 170 deg West (over the Pacific Ocean). The depot is assumed to consume its own monopropellant at a rate of 0.14kg/day for station keeping [35]. The manufacturing and operating costs of the depot are assumed to be $200M and $13,000/day, respectively.

The launch vehicle used for this analysis is based on a Falcon 9 launcher with an assumed maximum payload capacity of 8,300kg. A launcher is assumed to be available every 30 days for resupply of the depot. The mass-specific launch cost is assumed to be $11,300/kg.

Six different servicer designs are simulated in this paper based on their propulsion systems and the number of tools they are integrated with. We define three types of servicers based on the propulsion system they use to fly: high-thrust servicers; low-thrust servicers; and multimodal servicers, which integrate both high- and low-thrust propulsion technologies. We then define two subtypes based on the number of tools they are equipped with: 1 versatile servicer, and 4 specialized servicers. The versatile servicer can provide all seven defined services whereas the specialized servicers can only provide the services for which their tools are suited. The detailed assumptions are given in Tables 4, 5, and 6 for the high-thrust servicers, low-thrust servicers, and multimodal servicers, respectively. The baseline dry mass for the high-thrust versatile, low-thrust versatile, and multimodal specialized servicers is taken from Ref. [13]. This baseline mass is decreased for the high-thrust and low-thrust specialized servicers which are assumed to be less capable and smaller in size than their versatile counterparts. Similarly, this baseline mass is increased for the multimodal versatile servicers because they integrate an additional propulsion system compared to their high- and low-thrust versatile counterparts. The model for the low-thrust propulsion system is based on Northrop Grumman's MEV system of four XR-5 Hall Thrusters [6, 33].

Finally, this paper assumes the refueling of the servicers to be instantaneous operations. The justification behind this assumption is that, as OOS operations become routine, refueling of the servicers will likely not take more than one time step in the dynamic network (*i.e.*, 2 days) [36]. This assumption can be modified depending on the technology performance.

**Table 4. Assumptions related to the high-thrust servicers (*).**

| **Parameters** | **Versatile** | **Specialized 1** | **Specialized 2** | **Specialized 3** | **Specialized 4** |
|---|---|---|---|---|---|
| **Tools** | T1, T2, T3, T4 | T1 | T2 | T3 | T4 |
| **Dry mass [kg]** | 3,000 [13] | 2,000 | 2,000 | 2,000 | 2,000 |
| **Propellant capacity [kg]** | 1,000 | 1,000 | 1,000 | 1,000 | 1,000 |
| **Manufacturing cost [$M]** | 75 | 50 [37] | 50 [37] | 50 [37] | 50 [37] |
| **Operating cost [$/day]** | 13,000 [20] | 13,000 [20] | 13,000 [20] | 13,000 [20] | 13,000 [20] |
| **Propellant type** | Bi-propellant | Bi-propellant | Bi-propellant | Bi-propellant | Bi-propellant |
| **Specific Impulse [s]** | 316 | 316 | 316 | 316 | 316 |
| **Flight durations [days]** | 2, 4 | 2, 4 | 2, 4 | 2, 4 | 2, 4 |

* References are indicated within brackets. Data without reference is assumed.

**Table 5. Assumptions related to the low-thrust servicers (*).**

| | **Versatile** | **Specialized 1** | **Specialized 2** | **Specialized 3** | **Specialized 4** |
|---|---|---|---|---|---|
| **Tools** | T1, T2, T3, T4 | T1 | T2 | T3 | T4 |
| **Dry mass [kg]** | 3,000 [13] | 2,000 | 2,000 | 2,000 | 2,000 |
| **Propellant capacity [kg]** | 300 | 300 | 300 | 300 | 300 |
| **Manufacturing cost [$M]** | 75 | 50 [37] | 50 [37] | 50 [37] | 50 [37] |
| **Operating cost [$/day]** | 13,000 [20] | 13,000 [20] | 13,000 [20] | 13,000 [20] | 13,000 [20] |
| **Propellant type** | Low-thrust propellant | Low-thrust propellant | Low-thrust propellant | Low-thrust propellant | Low-thrust propellant |
| **Specific Impulse [s]** | 1,790 [6,33] | 1,790 [6,33] | 1,790 [6,33] | 1,790 [6,33] | 1,790 [6,33] |
| **Thrust [N]** | 1.74 [6,33] | 1.16 [6,33] | 1.16 [6,33] | 1.16 [6,33] | 1.16 [6,33] |
| **Flight durations [days]** | 10, 14, 30, 34 | 10, 14, 30, 34 | 10, 14, 30, 34 | 10, 14, 30, 34 | 10, 14, 30, 34 |

* References are indicated within brackets. Data without reference is assumed.

**Table 6. Assumptions related to the multimodal servicers (*).**

| | **Versatile** | **Specialized 1** | **Specialized 2** | **Specialized 3** | **Specialized 4** |
|---|---|---|---|---|---|
| **Tools** | T1, T2, T3, T4 | T1 | T2 | T3 | T4 |
| **Dry mass [kg]** | 4,000 | 3,000 [13] | 3,000 [13] | 3,000 [13] | 3,000 [13] |
| **Bi-propellant capacity [kg]** | 1000 | 1000 | 1000 | 1000 | 1000 |
| **Low-thrust propellant capacity [kg]** | 300 | 300 | 300 | 300 | 300 |
| **Manufacturing cost [$M]** | 100 | 75 | 75 | 75 | 75 |
| **Operating cost [$/day]** | 13,000 [20] | 13,000 [20] | 13,000 [20] | 13,000 [20] | 13,000 [20] |
| *High-thrust propulsion parameters* | | | | | |
| **Propellant type** | Bi-propellant | Bi-propellant | Bi-propellant | Bi-propellant | Bi-propellant |
| **Specific Impulse [s]** | 316 | 316 | 316 | 316 | 316 |
| **Flight durations [days]** | 2, 4 | 2, 4 | 2, 4 | 2, 4 | 2, 4 |
| *Low-thrust propulsion parameters* | | | | | |
| **Propellant type** | Low-thrust propellant | Low-thrust propellant | Low-thrust propellant | Low-thrust propellant | Low-thrust propellant |
| **Specific Impulse [s]** | 1,790 [6,33] | 1,790 [6,33] | 1,790 [6,33] | 1,790 [6,33] | 1,790 [6,33] |
| **Thrust [N]** | 1.74 [6,33] | 1.16 [6,33] | 1.16 [6,33] | 1.16 [6,33] | 1.16 [6,33] |
| **Flight durations [days]** | 10, 14, 30, 34 | 10, 14, 30, 34 | 10, 14, 30, 34 | 10, 14, 30, 34 | 10, 14, 30, 34 |

* References are indicated within brackets. Data without reference is assumed.

### iii. Case Studies' Scenarios

Using the assumptions presented previously, several scenarios are designed to demonstrate the framework's value in supporting trade studies related to the propulsion technologies of the servicers. Two case studies are considered in this paper: the short-term operational scheduling of a multimodal servicer; and the long-term strategic planning of six different OOS architectures. Three different market conditions are defined by considering 30 satellites (low demand rate), 71 satellites (medium demand rate), or 142 satellites (high demand rate).

The first case study aims to demonstrate the optimizer's ability to automatically trade between the high- and low-thrust engines of a multimodal servicer to maximize the short-term profits of an OOS venture. The regular scheduling of the short-term operations of OOS infrastructures will be essential to account for random demand (*e.g.*, repair need) and remain competitive. In this first case study, the servicer is assumed to be initially pre-deployed at some customer satellite to account for the fact that it may already be providing services when the re-planning event occurs. The framework is run for a single planning horizon of the RH procedure. This scenario is summarized in Table 7.

**Table 7. Scenario definition for the short-term operational scheduling case study.**

| | |
|---|---|
| **Servicers** | 1 multimodal versatile servicer |
| **Depot** | 1 depot pre-deployed at 170deg West on the GEO orbit |
| **Planning horizon** | 90 days |
| **Customer fleet** | 142 customer satellites (high demand) |

| **Initial conditions** | The servicer is initially deployed at the *SBIRS GEO 2* satellite due to a prior service. |
|---|---|

The second case study aims to compare the long-term value (*i.e.*, profits minus initial investments) of six different OOS architectures leveraging various propulsion technologies. This is done by running the OOS logistics framework for a 5-year time horizon, while leveraging the RH procedure embedded in the framework to account for the random service needs. For this case study, we assess the value of six different OOS architectures under low service demand (30 satellites), medium service demand (71 satellites), and high service demand (142 satellites). The architectures, given in Table 8, are defined based on the servicers' propulsion technologies (*i.e.*, high-thrust, low-thrust, or multimodal propulsion) and platform designs (*i.e.*, versatile or specialized). Architectures involving 1 versatile servicer are referred to as *monolithic*, while those involving 4 specialized servicers are referred to as *distributed*. As in the first case study, an orbital depot is pre-deployed at a longitude of 170deg West on the GEO orbit. Finally, note that, although this is not demonstrated in this paper, the proposed framework can simulate OOS architectures involving servicers with different propulsion technologies (*e.g.*, low-thrust and high-thrust servicers working in concert to provide services).

**Table 8. Definition of the OOS architectures for the long-term strategic planning case study.**

| **Architecture parameters** | **Architecture 1 (high-thrust monolithic)** | **Architecture 2 (high-thrust distributed)** | **Architecture 3 (low-thrust monolithic)** | **Architecture 4 (low-thrust distributed)** | **Architecture 5 (multimodal monolithic)** | **Architecture 6 (multimodal distributed)** |
|---|---|---|---|---|---|---|
| **Number of Servicers** | 1 | 4 | 1 | 4 | 1 | 4 |
| **Servicer Propulsion** | **H**igh-thrust | **H**igh-Thrust | **L**ow-thrust | **L**ow-Thrust | **M**ultimodal | **M**ultimodal |
| **Servicer Platform Design** | **V**ersatile (4 tools per servicer) | **S**pecialized (1 tool per servicer) | **V**ersatile (4 tools per servicer) | **S**pecialized (1 tool per servicer) | **V**ersatile (4 tools per servicer) | **S**pecialized (1 tool per servicer) |

### B. Case Study 1: Short-Term Operational Scheduling

Whenever an OOS operator decides to re-plan the operations of their infrastructure, they would run the proposed framework with an updated set of service needs and the initial conditions corresponding to the state of the infrastructure at re-planning. In the output, they would then get a breakdown of the servicers' operations over the planning horizon to inform short-term decision making.

The present case study is run for a customer base of 142 customer satellites. However, not all of these satellites will need a service during the period considered for re-planning. The framework thus automatically includes in the network only the satellites that display service needs over the planning horizon. This minimizes the size of the optimization problem thus allowing for a shorter computational time. For example, in the present case study, out of the 142 customer satellites, only 19 exhibit at least one service need, 5 of which are actually visited by the multimodal versatile servicer in the optimal solution. Figure 7 gives the network for the considered scenario. The customer satellites are represented by colored dots. The red dots correspond to the visited customer satellites, whereas the blue dots are not visited by the servicer. The black star in Fig. 7 represents the orbital depot.

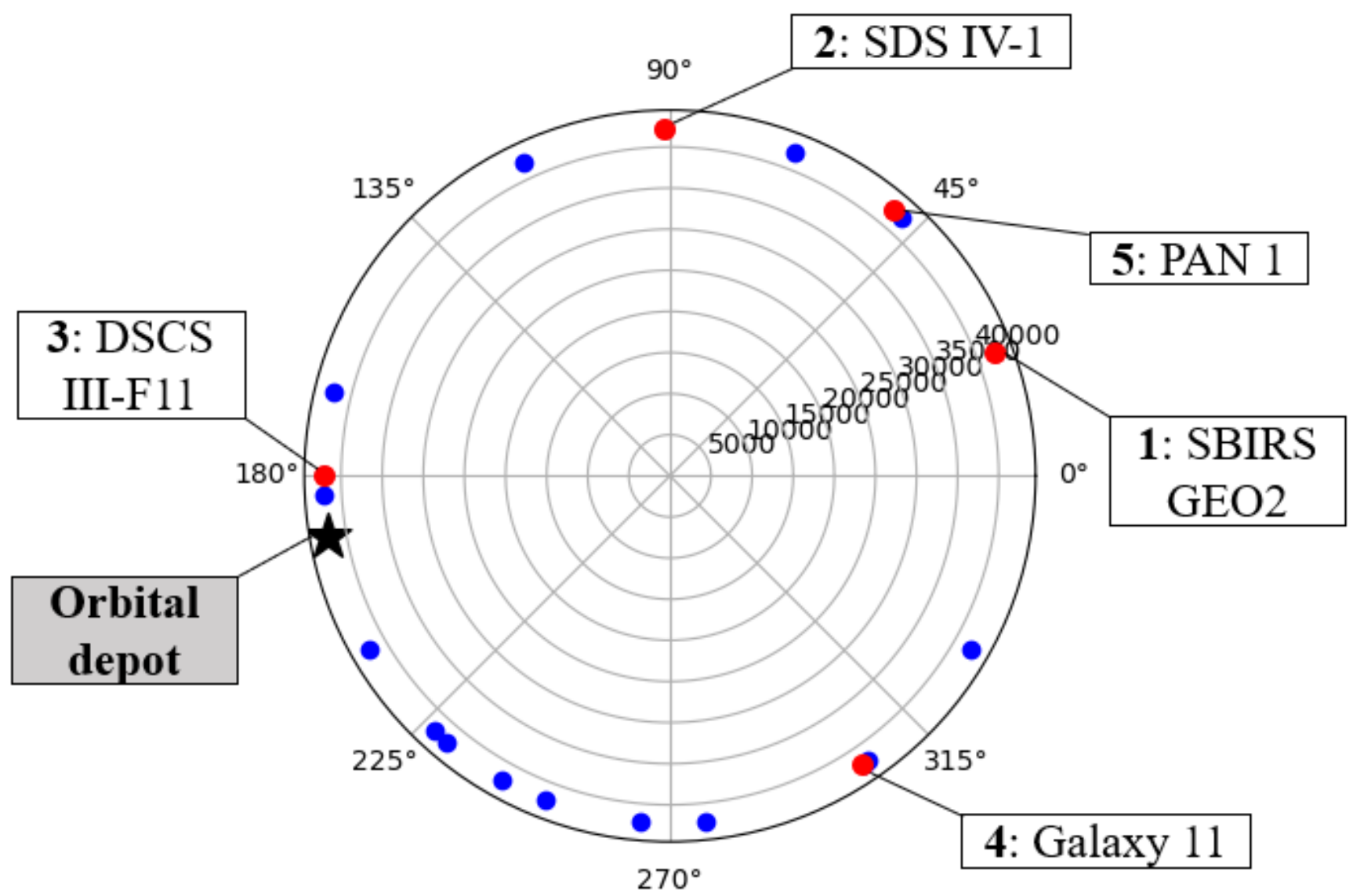

**Fig. 7 Static network used by the optimizer to plan the short-term operations of a multimodal versatile servicer.**

The optimization framework is run on an Intel® Core™ i7-9700, 3.00GHz platform with the Gurobi 9 optimizer. The solution was reached in 45 seconds with a gap of 1% a stopping criterion. Figures 8-10 illustrate the optimal operations of the multimodal versatile servicer over time.

As seen in Fig. 8-10, the optimizer has the servicer fly with both high- and low-thrust propulsion systems depending on the time the servicer has left to travel between the nodes of the network. The high- and low-thrust transportation arcs of the servicer are represented with different dashed lines on the figures. Note that at $t = 50$, the servicer leaves the *DSCS III-F11* satellite for the depot, and gets refueled in both low-thrust propellant and bi-propellant in order to perform the next low-thrust maneuver at $t = 54$ and high-thrust maneuver at $t = 70$.

For information, the operations presented in Fig. 8-10 lead to a maximized profit of $58.8M over this 90-day period. The servicer was able to capture $75M in revenues from the provision of 4 services (to, in chronological order, *SDS IV-1*, *DSCS III-F11*, *Galaxy 11*, and *PAN 1*) for a total cost of $16.2M, including launch, purchase of propellant, delay penalties, and infrastructure operations. Note that these values may not be final for this 90-day period as random service needs may occur within that period, thus triggering future re-planning events and likely modifying the short-term profits.

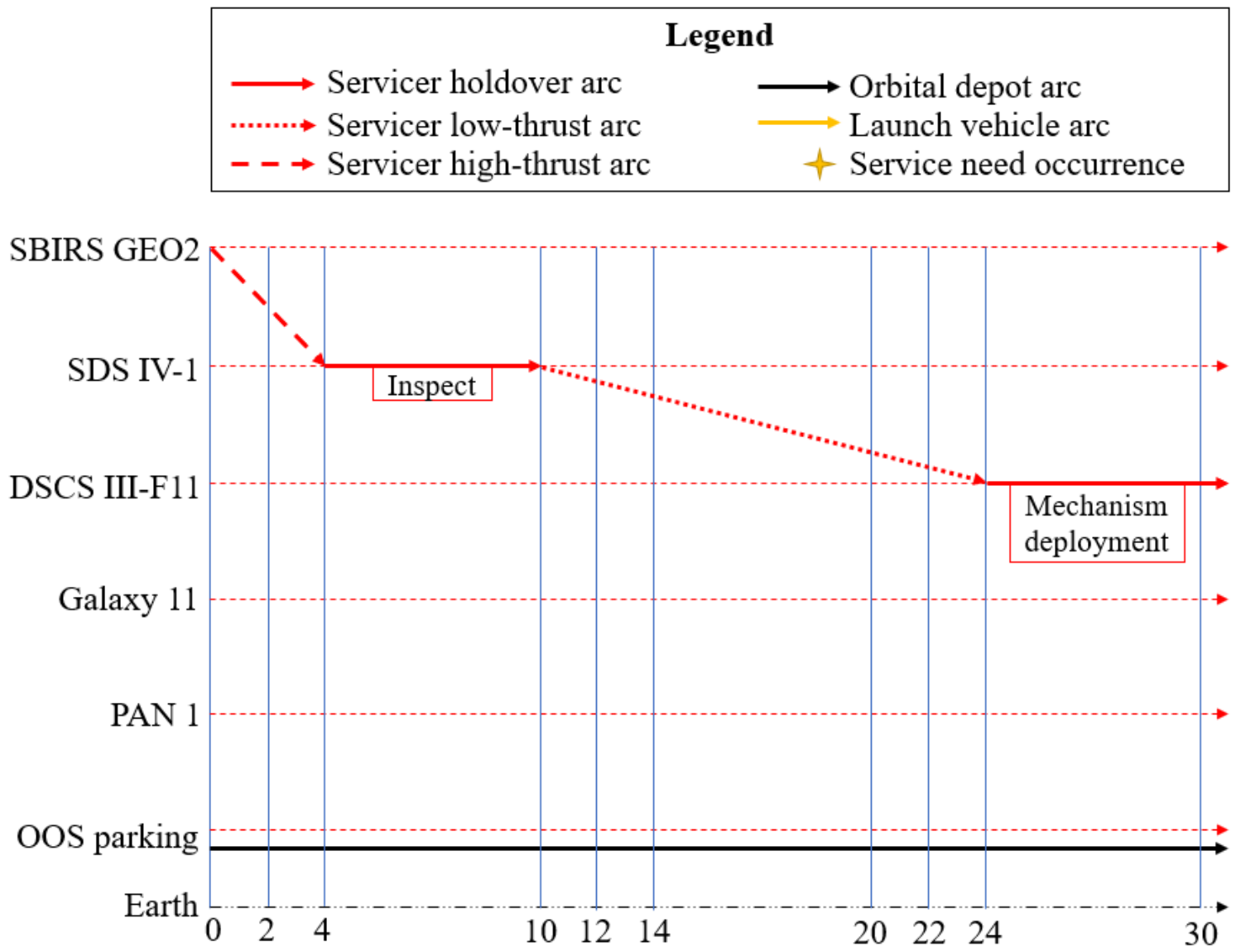

**Fig. 8 Breakdown over the first month (days 0-30) of the operations of an OOS infrastructure consisting of a multimodal versatile servicer and an orbital depot.**

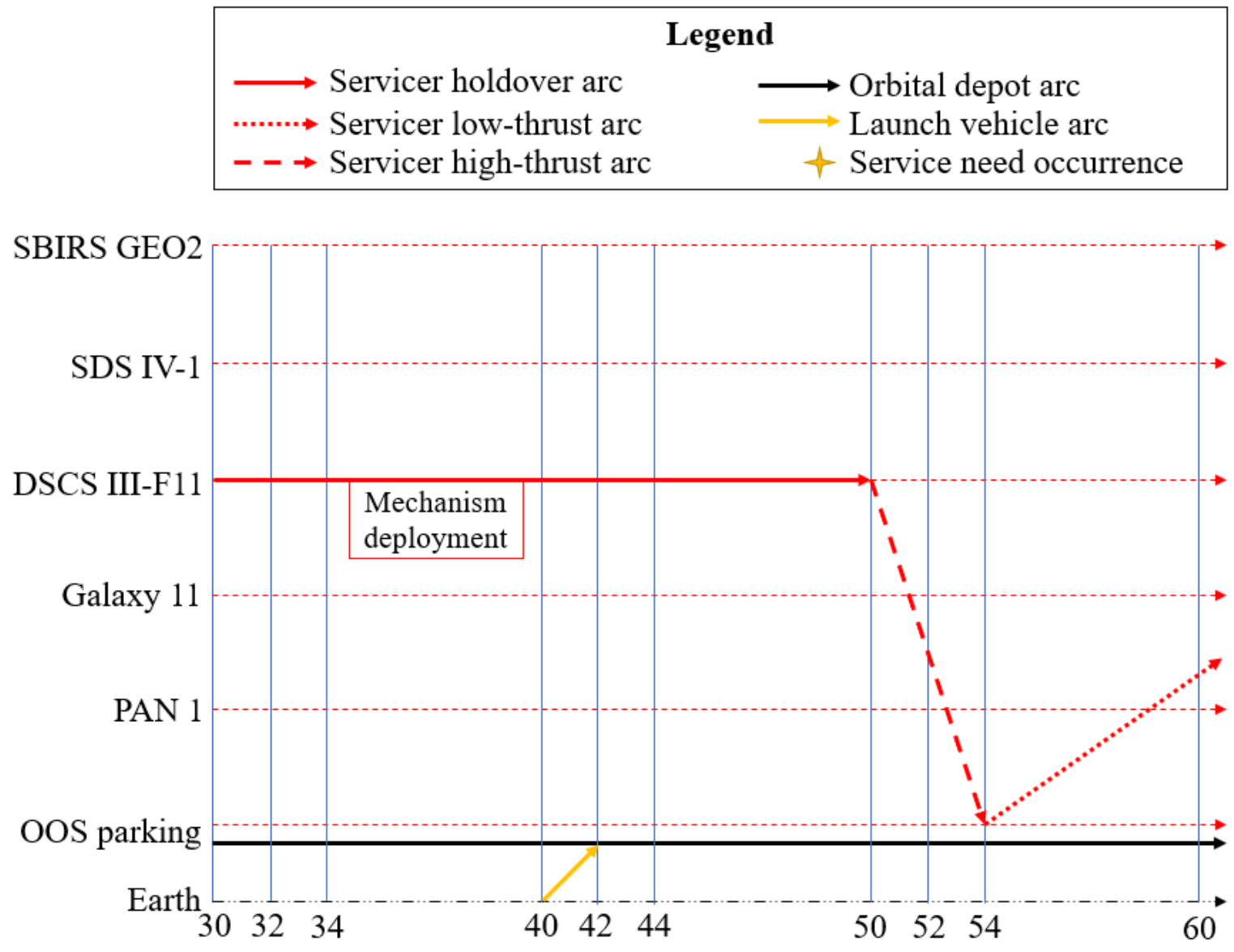

**Fig. 9 Breakdown over the second month (days 30-60) of the operations of an OOS infrastructure consisting of a multimodal versatile servicer and an orbital depot.**

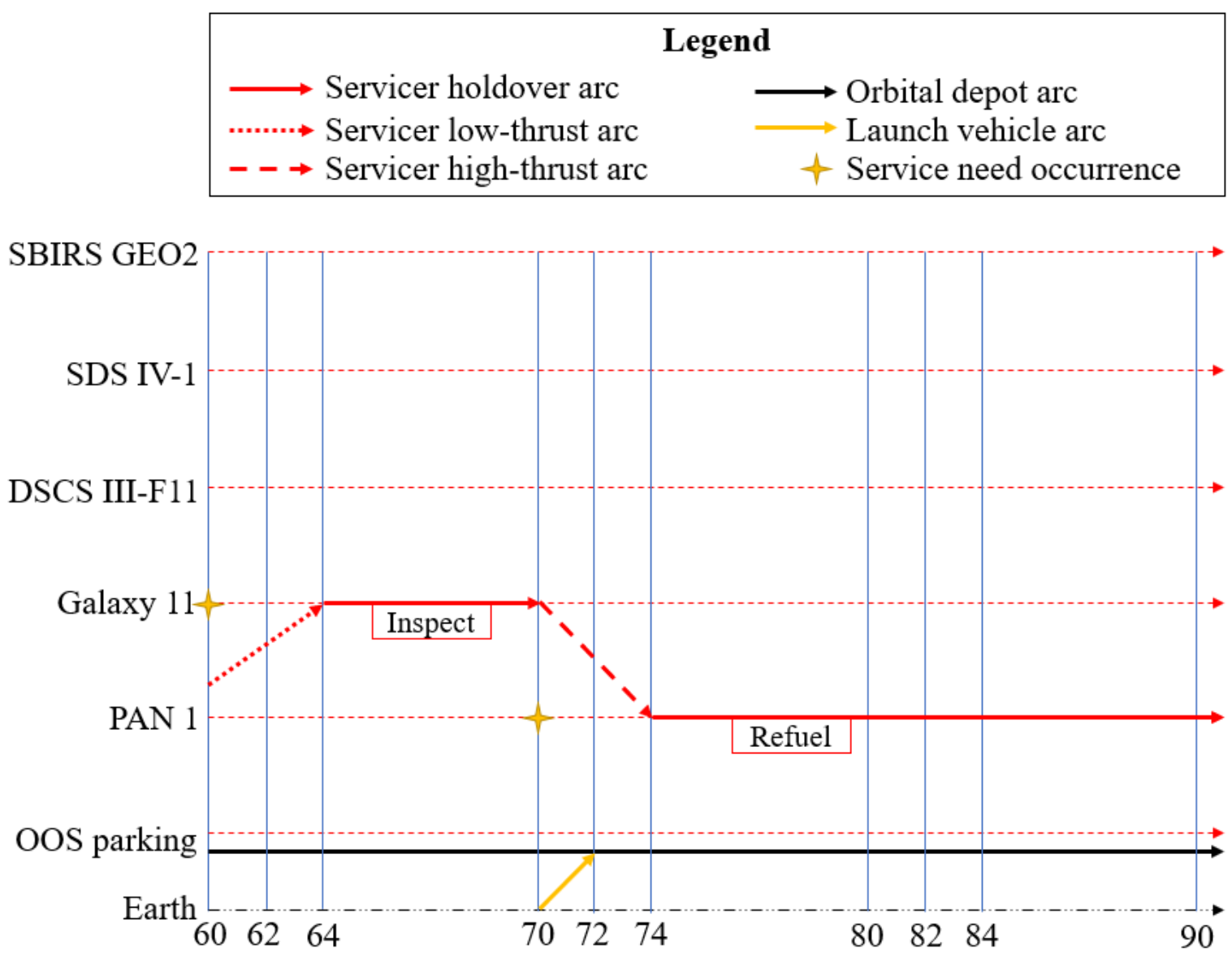

**Fig. 10 Breakdown over the third month (days 60-90) of the operations of an OOS infrastructure consisting of a multimodal versatile servicer and an orbital depot.**

### C. Case Study 2: Long-Term Strategic Planning

Careful long-term strategic planning is critical for new entrants in the OOS market to design servicing infrastructures that are resilient to competition and variations in service demand, and are capable of generating significant marginal profits in the prevision of market growth. Decision-makers will need to explore a large tradespace capturing both the operational strategies and the designs of the servicers and orbital depots, while considering the uncertainties in service demand. The generalization of the previous framework [22] enabling the tradeoff between high- and low-thrust technologies is believed to be a significant step forward toward that goal. On the other hand, the RH procedure is critical to properly account for random service needs, which, in essence, are nothing more than unplanned sources of revenues.

For this case study, three different analyses are performed. First, we compare the values of the monolithic and distributed architectures for each propulsion technology and different market conditions (30, 71, 142 customer satellites). Then, we discuss the impact of the servicers' propulsion technologies on the value of the OOS infrastructures. Finally, we take a look at how sensitive the value of an OOS infrastructure can be to variations in the mass of the servicers.

#### i. Trading Architectural Options: Monolithic vs. Distributed Architectures

The purpose of this section is to compare the performance of the monolithic and distributed architectures per propulsion system. More specifically, we compare Architecture 1 with Architecture 2, Architecture 3 with Architecture 4, and Architecture 5 with Architecture 6 (*cf* Table 8) for three different levels of service demands (30, 71, 142 customer satellites). For this analysis, 18 simulations are run on an Intel® Core™ i7-9700, 3.00GHz platform with the Gurobi 9 optimizer. The simulation times range between 1 and 129 minutes with an average MILP gap between 0 and 0.13%.

The results for the values of Architecture 1 (*i.*e., 1 high-thrust versatile servicer) and Architecture 2 (*i.e.*, 4 high-thrust specialized servicers) are first presented in Fig. 11 for the three different levels of service demands. As can be seen on this figure, the monolithic architecture is more valuable than the distributed one over the 5-year time horizon for a small OOS market (*e.g.*, 30 satellites). This is due to a higher initial investment and operating cost of the distributed architecture. With a small market, the OOS infrastructures cannot pay back their initial investments in less than 5 years of business operations. However, as the market grows to a customer base of 71 satellites, the distributed architecture is seen to be catching up with the value of the monolithic architecture. Under this market condition, the monolithic architecture is barely valuable after 5 years of operations while the distributed architecture is still $100M away from the breakeven. Finally, for a large market of 142 customer satellites, both architectures become valuable after 2.5 years of operations with the distributed architecture being more profitable than the monolithic one. This is because, unlike the monolithic architecture, the distributed architecture can *parallelize* the provision of services with 4 different servicers thus allowing for a larger profitability rate. Essentially, this analysis shows that the distributed architecture has a larger potential in marginal profitability as the market grows than the monolithic architecture. Note that the results presented in Fig. 11 slightly differ from those presented in Ref. [22] because the high-thrust servicers are given the option to fly over 4-day trajectories instead of over 2-day trajectories only.

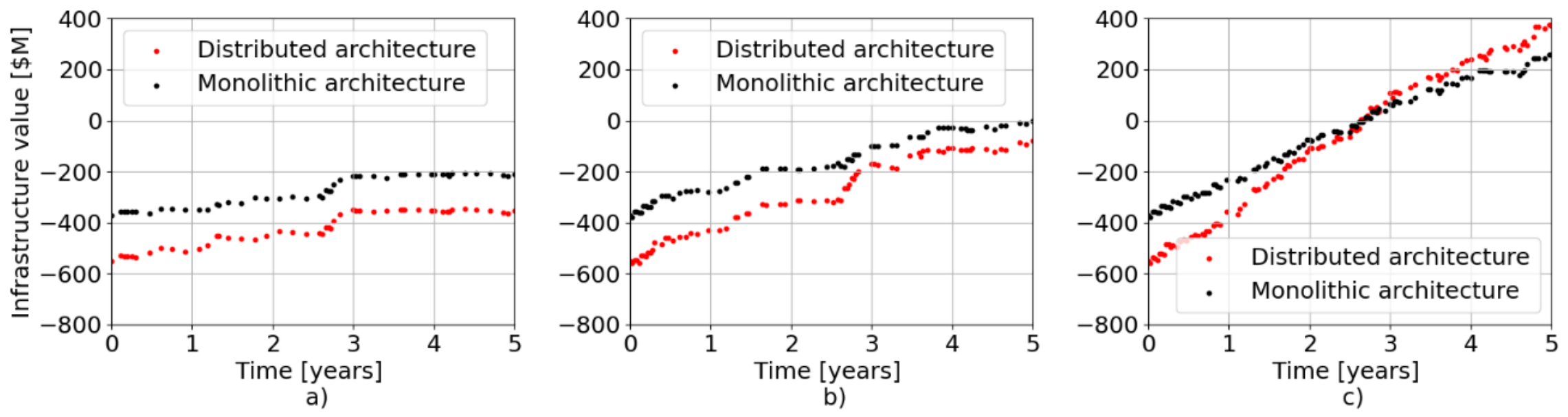

**Fig. 11 Values of the high-thrust monolithic and distributed architectures (Architecture 1 vs Architecture 2) for: a) 30 satellites; b) 71 satellites; c) 142 satellites**

A similar analysis is carried out to compare the values of Architecture 3 (*i.e.*, 1 low-thrust versatile servicer) and Architecture 4 (*i.e.*, 4 low-thrust specialized servicers). The results are presented in Fig. 12 for the three different levels of service demands. A similar trend as with the high-thrust architectures can be observed on those figures. However, the distributed architecture catches up faster with the value of the monolithic architecture as the market grows. This can be seen in Fig. 12 with a customer base of 71 satellites: the low-thrust distributed architecture becomes more valuable than its monolithic counterpart after 2.8 years of operations. As justified for the high-thrust architectures, the low-thrust distributed architecture has a greater potential in marginal profitability as the market grows due to the parallelization of the services.

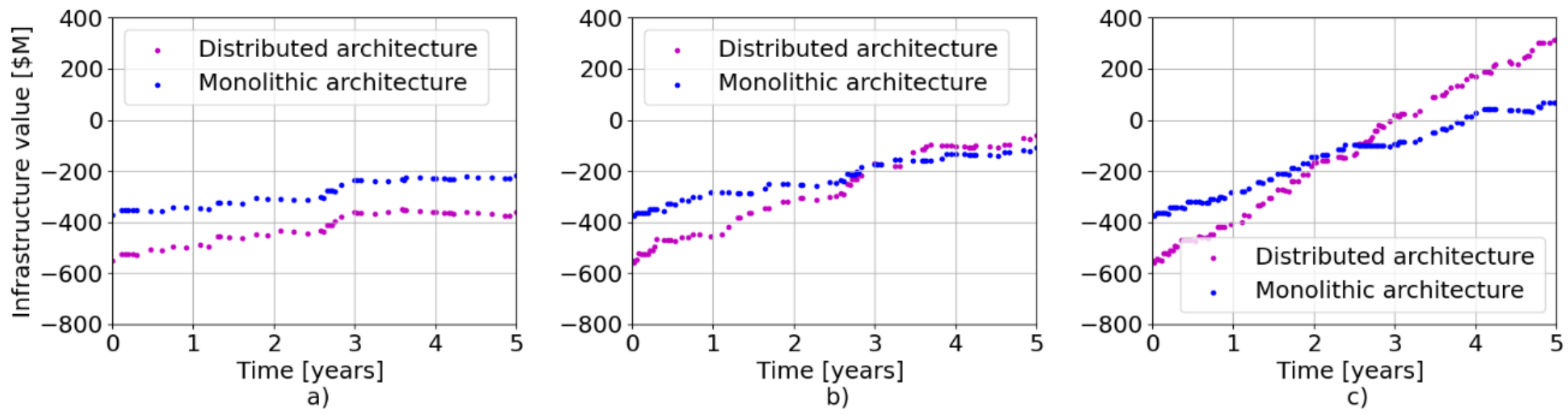

**Fig. 12 Values of the low-thrust monolithic and distributed architectures (Architecture 3 vs Architecture 4) for: a) 30 satellites; b) 71 satellites; c) 142 satellites**

Finally, we compare Architecture 5 (*i.e.*, 1 multimodal versatile servicer) and Architecture 6 (*i.e.*, 4 multimodal specialized servicers). The results are presented in Fig. 13 for the three different levels of service demands. One can observe the same trends as with the high- and low-thrust architectures. Again, the multimodal distributed architecture has a greater potential in marginal profitability as the market grows due to the parallelization of the services.

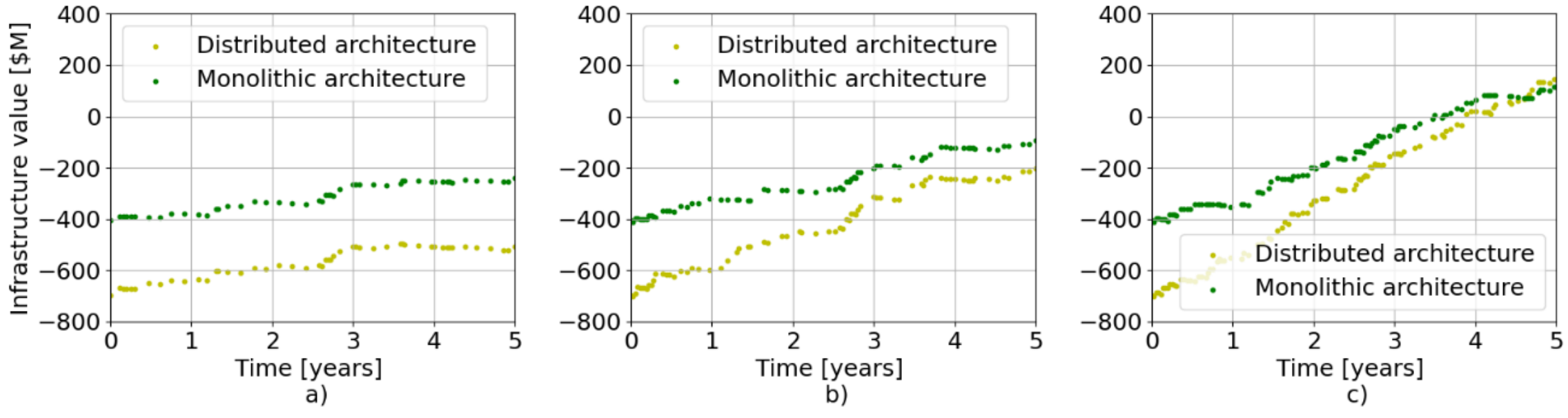

**Fig. 13 Values of the multimodal monolithic and distributed architectures (Architecture 5 vs Architecture 6) for: a) 30 satellites; b) 71 satellites; c) 142 satellites**

### ii. Trading Propulsion Options: High-Thrust vs. Low-Thrust vs. Multimodal Propulsion

This second analysis aims to demonstrate through an example how the proposed OOS logistics framework can be leveraged to inform decision making about the propulsion systems of the servicers. The results presented here are for an assumed large market of 142 customer satellites, and are sourced from the 18 simulations introduced in IV.C.i. Note that the analysis which follows is valid for the assumptions made in subsection IV.A. Different assumptions can easily be input to our framework, leading to potentially different conclusions from those presented below.

The focus for this example analysis is on whether using multimodal servicers significantly contributes to the value of the OOS infrastructure. Multimodal servicers offer more flexibility to OOS companies through an optimal balance between high-thrust propulsion (*i.e.*, better responsiveness but higher operating cost) and low-thrust propulsion (*i.e.*, slower but cheaper operations). However, this comes at a cost with respect to the high- and low-thrust servicers: a higher initial investment (+$25M), and an added dry mass (+1,000kg), which in turn either increases the operating costs or reduces the revenues available to the servicers.

To give some insight into this important tradeoff, we first take a look at the propulsion system utilization ratio of the multimodal versatile and specialized servicers presented in Fig. 14. This metric gives the OOS planners insight into how frequently the multimodal servicers use their high-thrust and low-thrust engines. In Fig. 14, one can see that the multimodal specialized servicers (*i.e.*, Architecture 6) use their low-thrust engines most of the time. Indeed, the parallelization of services enabled through the use of multiple servicers gives them more time to perform the rendezvous maneuvers. The optimizer thus tends to favor the low-thrust engine over the high-thrust engine due to a lower propellant consumption. Based on this observation, OOS planners may decide to deploy low-thrust specialized servicers instead of multimodal specialized servicers to simplify and reduce the cost of manufacturing and operations.

Now take a look into the propulsion system utilization ratio for the multimodal versatile servicer (Architecture 5). As can be seen in Fig. 14, the servicer uses its high-thrust engine almost as much as its low-thrust engine, which proves the two technologies are similarly valuable for this type of OOS architecture employing a single versatile servicer. Indeed, in this case, the parallelization of services is inexistent since there is only one servicer to address all service needs in a sequential manner. The optimizer thus makes the most of both high- and low-thrust engines to find the optimal balance between servicing responsiveness and operating costs: the more responsive the infrastructure, the more revenues it generates but also the more costly it becomes to operate. Clearly, the tradeoff between the multimodal versatile servicer and its high- and low-thrust counterparts is not as straightforward as in the case of multimodal specialized servicers, and Fig. 14 alone is not sufficient to properly inform decision making.

To get a better understanding of the tradeoff between the propulsion technologies of the versatile servicer, we need to consider this problem from a broader perspective. To that end, Fig. 15 gives the cost breakdown for Architecture 1 (*i.e.*, 1 high-thrust versatile servicer), Architecture 3 (*i.e.*, 1 low-thrust versatile servicer), and Architecture 5 (*i.e.*, 1 multimodal versatile servicer). As a reminder, the *value* metric is defined as the *revenues* minus the sum of the *initial investments* and *operational costs*. As can be seen in Fig. 15, employing a high-thrust versatile servicer is significantly more valuable for a large OOS market (*i.e.*, 142 customer satellites) than the other two options. The multimodal versatile servicer is $143M short of the high-thrust servicer in terms of value. This seems counterintuitive as the multimodal servicer can benefit from its low-thrust engine. However, it also is 1,000kg heavier, by assumption, than the high-thrust servicer due to the additional propulsion system. A heavier multimodal servicer leads to a larger bi-propellant consumption and may not be able to fly along low-thrust transportation arcs with relatively low *servicer mass upper bounds* (*cf* section II.B and Eq. 26). Architecture 3 (*i.e.*, 1 low-thrust servicer) is the least valuable of the three options investigated here due to low servicing responsiveness and thus lower revenues overall.

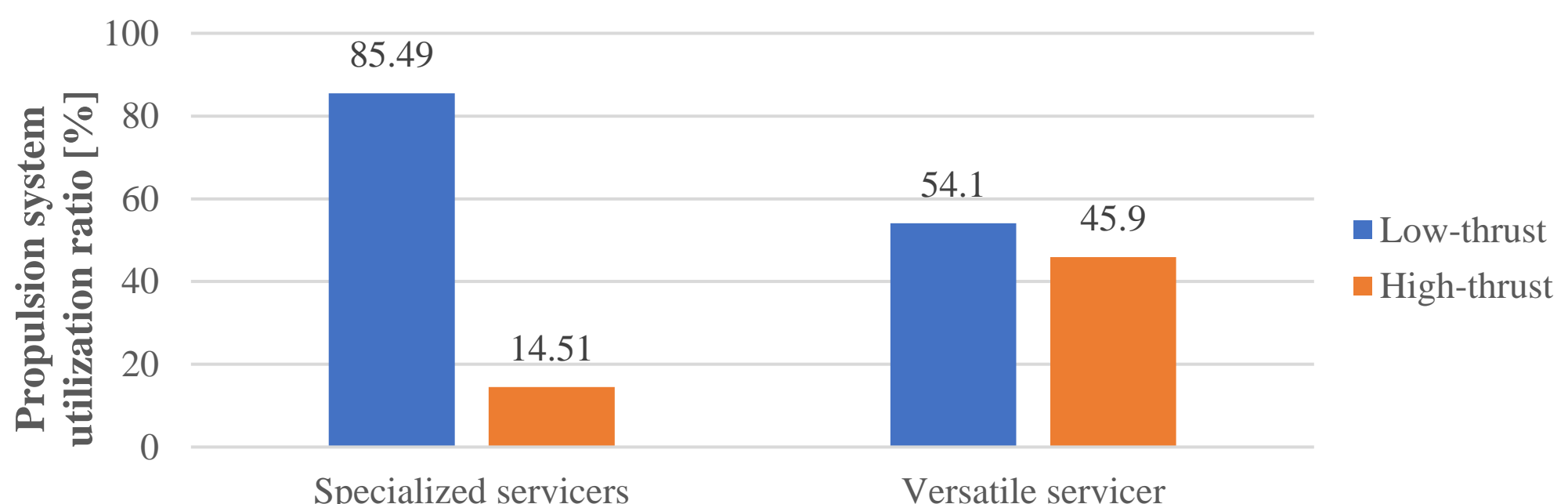

**Fig. 14 Propulsion system utilization ratio of the multimodal versatile servicer (Architecture 5) and multimodal specialized servicers (Architecture 6).**

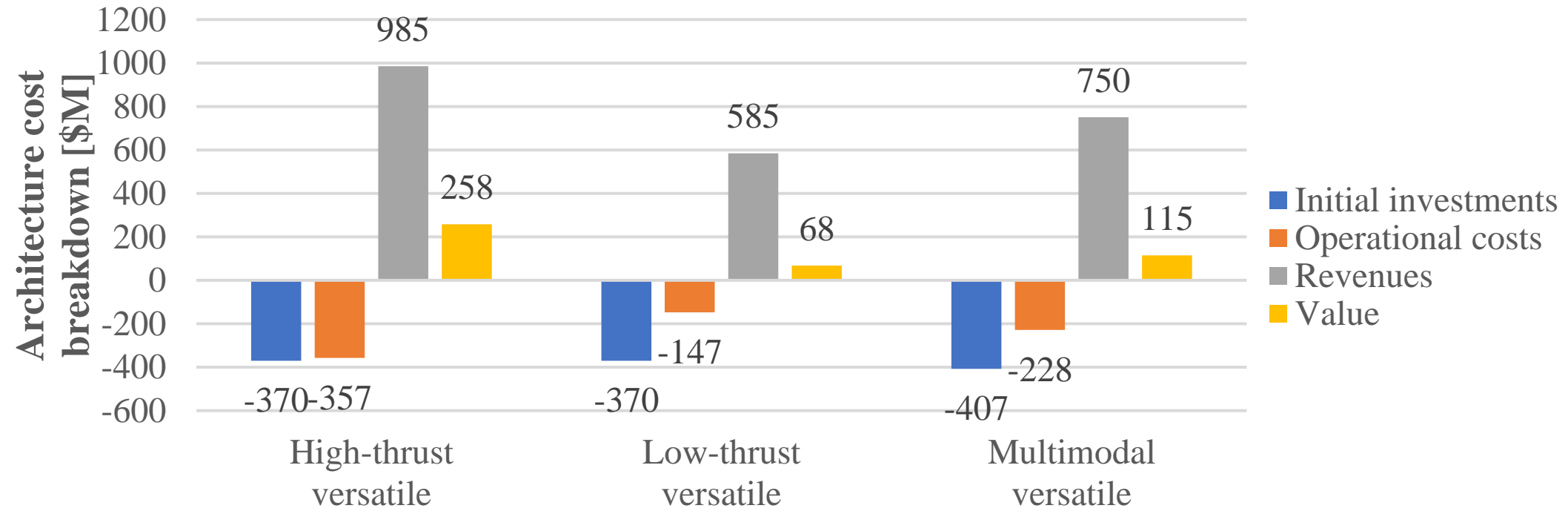

**Fig. 15 Cost breakdown for Architecture 1 (high-thrust versatile servicer), Architecture 3 (low-thrust versatile servicer), and Architecture 5 (multimodal versatile servicer) after 5 years of operations.**

### iii. Sensitivity Analysis with Respect to Servicer Mass

In this subsection, we present how sensitive the revenues, operating costs, and profits are to variations in the mass of the Power Processing Unit (PPU) of the multimodal versatile servicer in Architecture 5. In this analysis, the initial investments to manufacture and deploy the OOS infrastructure into space are not included in the results to fairly assess the profit-generating performance of the servicer as a function of its dry mass. The long-term strategic planning of Architecture 5 is run for a servicer mass of 3,000kg, 4,000kg, 5,000kg, 6,000kg, an OOS market of 142 customer satellites, and a 5-year timeline. Each simulation ran in less than 40 minutes on an Intel® Core™ i7-9700, 3.00GHz platform with the Gurobi 9 optimizer and a MILP gap of 1% a stopping criterion.

Figure 16 gives the revenues, the operating costs, and the profits of the architecture which are plotted with the same Y-axis scale to facilitate the comparison between their rates of change. The analysis which follows is valid for the assumptions made in this paper and should not be used to draw conclusions about the general design of sustainable OOS infrastructures. Instead, OOS planners would need to input their own trustworthy data to the proposed framework, which is highly and easily parameterizable.

The general trend observed on those figures is that decreasing the mass of a servicer has a larger impact on the revenue stream and profits than on the operating costs. For example, after 5 years of operations, the 3000kg servicer has generated about $80 million in excess of those generated by the other 3 servicers, while the operating costs are similar for all four servicers. Indeed, as the mass of a servicer decreases, the number of feasible rendezvous maneuvers increases as modeled by Eq. 26 and the concept of *servicer mass upper bound*. Decreasing the mass of a servicer thus unlocks additional sources of revenues available to the OOS infrastructures by enabling additional feasible trajectories. On the other hand, the revenues generated by the 4,000kg-to-6,000kg servicers are almost identical, which shows that there is a servicer mass threshold between 3,000kg and 4,000kg beyond which the subset of satellites accessible to the servicer through feasible maneuvers remains principally unchanged. In other words, the 3,000kg servicer taps into more revenue sources than the 4,000kg-to-6,000kg servicers.

Interestingly, the operating costs are almost identical no matter the mass of the servicer. Indeed, a lighter servicer will provide more services and thus will need to fly more often than a heavier servicer. However, flying a lighter servicer will require less propellant overall per flight than a heavier servicer.

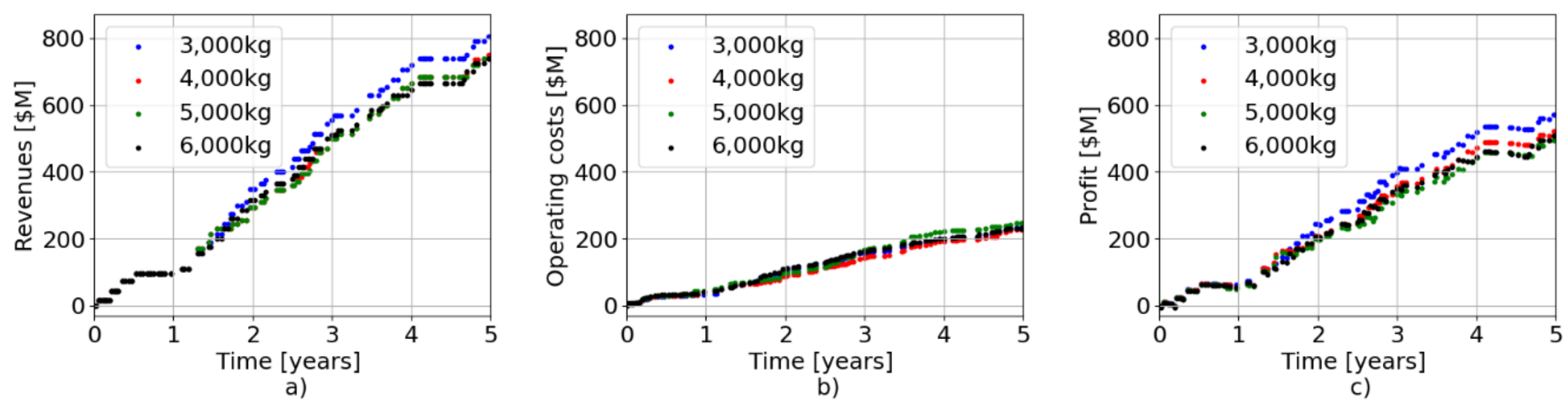

**Fig. 16 a) Revenues; b) Operating costs; and c) profit of Architecture 5 for four different masses of the multimodal versatile servicer.**

This analysis gives some insight into how valuable incremental mass-reduction can be to an OOS company. The profits generated by reducing the mass of a servicer could then be reinvested into the company to capture larger shares of the OOS market.

### iv. Discussions

The analyses presented in this paper do not aim to provide definitive answers regarding the design and operations of OOS infrastructures. Instead, they intend to show the generality of the proposed approach and its ability to perform an automatic tradeoff between the high- and low-thrust propulsion technologies. This method could be used to investigate the long-term value of simpler near-term infrastructures, and could easily be modified to assess when and how a given OOS infrastructure should be upgraded in response to a dynamically evolving OOS market.

These analyses are typical examples of how OOS designers could use the proposed framework to effectively inform decision management for the long term. This can be done by running sensitivity analyses over various servicer and depot designs, operational schemes, and even over the nature and price tags of the commodities needed to support the operations of OOS infrastructures. For example, how would the results presented in this section change if a depot is not deployed at all, or its resupply by the launch vehicle is not as frequent as assumed in this paper (*i.e.*, 30 days)? Wouldn't low-thrust servicers become more valuable than their high-thrust and multimodal counterparts? Is there a situation where the multimodal versatile servicers are more advantageous than high-thrust versatile ones? And what if, as lunar water extraction becomes mainstream, servicers are designed and deployed with water-based low-thrust engines? Xenon gas is an expensive commodity, but water extracted from the Moon may bring in a whole new perspective on the design and operations of OOS infrastructures. As the cislunar economy develops and technology options diversify, OOS planners will need a tool to take a wide variety of parameters into account and come up with competitive and efficient business models. We believe the OOS logistics framework proposed in this paper is a good candidate tool for use by the industry.

## V. Conclusion

This paper generalizes the state-of-the-art OOS logistics method by incorporating the automatic tradeoff between the high- and low-thrust propulsion systems of the servicers. The proposed OOS logistics framework is capable of modeling and simulating complex sustainable OOS infrastructures that involve all kinds of depot and servicer models. More specifically, this paper develops three generic models of servicers: high-thrust-only servicers, low-thrust-only servicers, and multimodal servicers, which embed both high- and low-thrust engines. In addition, we introduce the concept of trajectory plugins that allow OOS planners to build their own high- and low-thrust trajectory models, effortlessly interface them with the framework, and test them as alternative operational strategies for the servicers. The developed piecewise linear approximation model can convert the nonlinear mass transformation model into the MILP-based OOS logistics optimization. Once the trajectory plugins and servicer models are defined, the framework leverages Mixed-Integer Linear Programming and the Rolling Horizon approach to solve the optimal short-term operational scheduling and long-term strategic planning of OOS infrastructures under various states of the OOS market.

The proposed OOS logistics framework makes two fundamental assumptions: (1) the customer satellites are distributed over a common circular orbit; and (2) the orbital depots are staged along the same circular orbit as the fleet of customer satellites. Because GEO servicing falls well within these assumptions, the framework is used to run two different case studies related to this new industry. First, we optimize the operations of a multimodal versatile servicer and an orbital depot over a 90-day period, and show that the optimizer is capable of selecting for each rendezvous maneuver the best propulsion technology and trajectory to maximize the short-term profits. The second case study deals with the long-term strategic planning of 6 different OOS architectures under 3 different OOS market conditions. We show that the proposed method can effectively support long-term decision making regarding the architectural paradigm of OOS infrastructures (*e.g.*, monolithic vs. distributed) and the propulsion systems to adopt for the servicers.

We believe that the framework proposed in this paper constitutes a significant step toward effectively guiding the decisions related to the nascent but promising OOS industry. A natural extension of the current framework will relax the assumption that the customer satellites are located along the same circular orbits. This will allow OOS planners to explore the design of infrastructures dedicated to other OOS markets than the GEO market.

## Appendix: Definitions and Assumptions for Case Study

**Table A1. Definitions of the main parameters used to model the service needs [22].**

| Parameter | Applies to deterministic and/or random needs? | Definition |
|---|---|---|
| *Revenue* | Deterministic and random | Revenue received by the OOS infrastructure for providing a service in response to some service need. |
| *Delay penalty cost* | Deterministic and random | How much it costs daily to the OOS operator to delay the beginning date of the service. |
| *Service duration* | Deterministic and random | The time it takes to provide a service. |
| *Service window* | Deterministic and random | Time interval within which the service must start, provided that the OOS operator decides to provide the service. |
| *Interoccurrence time* | Deterministic only | Used to generate service needs regularly spaced in time. |
| *Mean interoccurrence time* | Random only | Used to generate service needs according to a Poisson probability distribution. |

**Table A2. Notional service-tool mapping (represented by an 'X' in the table) [22].**

| Service type | T1: Refueling apparatus | T2: Observation sensors | T3: Dexterous robotic arm | T4: Capture mechanism |
|---|---|---|---|---|
| *Inspection* | | X | | |
| *Refueling* | X | | | |
| *Station Keeping* | | | | X |
| *Repositioning* | | | | X |
| *Retirement* | | | | X |
| *Repair* | | | X | |
| *Mechanism Deployment* | | | X | |

## Acknowledgments

This work is supported by the Defense Advanced Research Project Agency Young Faculty Award D19AP00127. The content of this paper does not necessarily reflect the position or the policy of the U.S. Government, and no official endorsement should be inferred. This paper has been approved for public release; distribution is unlimited.

## References

[1]. Henry, C., Geostationary satellite orders bouncing back. SpaceNews. https://spacenews.com/geostationary-satellite-orders-bouncing-back/#:~:text=Villain%20said%20the%20average%20capital,than%20%24100%20million%20to%20deploy. [Accessed 1/26/2021]

[2]. OSAM-1: On-Orbit Servicing, Assembly and Manufacturing. NASA. https://nexis.gsfc.nasa.gov/OSAM-1.html [Accessed 1/26/2021].

[3]. Parrish, J., Robotic Servicing of Geosynchronous Satellites (RSGS). Defense Advanced Research Project Agency. https://www.darpa.mil/program/robotic-servicing-of-geosynchronous-satellites [Accessed 2/24/2020].

[4]. SpaceLogistics. Northrop Grumman. https://www.northropgrumman.com/space/space-logistics-services/ [Accessed 1/26/2021].

[5]. NSR press releases. Northern Sky Research. https://www.nsr.com/nsr-report-forecasts-4-5-billion-in-cumulative-revenues-from-in-orbit-satellite-services-by-2028/ [Accessed 1/26/2021]

[6]. Aerojet Rocketdyne propulsion helps enable new satellite servicing market. Aerojet Rocketdyne. https://www.rocket.com/media/news-features/ar-propulsion-helps-enable-new-satellite-servicing-market#:~:text=The%20XR%2D5%20thrusters%20provide,docked%20with%20the%20host%20satellite.&text=These%20thrusters%20provide%20full%20six%20degree%2Dfreedom%20control%20for%20docking%20maneuvers. [Accessed 1/26/2021]

[7]. Xenon Acquisition Strategies for High-Power Electric Propulsion NASA Missions. NASA. https://ntrs.nasa.gov/api/citations/20150023079/downloads/20150023079.pdf [Accessed 1/26/2021]

[8]. Convert prices and calculate cost of materials and compounds. Aqua-calc. https://www.aqua-calc.com/calculate/materials-price [Accessed 1/26/2021]

[9]. Werner, D., Orbit Fab to launch first fuel tanker in 2021 with Spaceflight. SpaceNews. https://spacenews.com/orbit-fab-to-launch-with-spaceflight/ [Accessed 1/31/2021]

[10]. Monomethylhydrazine. Astronautix.com. http://www.astronautix.com/m/mmh.html [Accessed 1/26/2021]

[11]. Zhao, S., Gurfil, P., and Zhang, J., “Optimal servicing of geostationary satellites considering earth’s triaxiality and lunisolar effects,” *Journal of Guidance, Control, and Dynamics*, Vol. 39, 2016, pp. 2219–2231. https://doi.org/10.2514/1.G001424

[12]. Han, C., and Wang, X., "On-orbit servicing of geosynchronous satellites based on low-thrust transfers considering perturbations," *Acta Astronautica*, vol. 159, pp. 658-675, 2019.
https://doi.org/10.1016/j.actaastro.2019.01.041

[13]. Verstraete, A. W., Anderson, D., St. Louis, N. M., and Hudson, J., “Geosynchronous Earth Orbit Robotic Servicer Mission Design,” *Journal of Spacecraft and Rockets*, Vol. 55, 2018, pp. 1444–1452.
https://doi.org/10.2514/1.A33945

[14]. Yao, W., Chen, X., Huang, Y., and van Tooren, M., “On-orbit servicing system assessment and optimization methods based on lifecycle simulation under mixed aleatory and epistemic uncertainties,” *Acta Astronautica*, Vol. 87, 2013, pp. 8-13.
https://doi.org/10.1016/j.actaastro.2013.01.012

[15]. Sarton du Jonchay, T., and Ho, K., “Quantification of the responsiveness of on-orbit servicing infrastructure for modularized earth orbiting platforms,” *Acta Astronautica*, Vol. 132, 2017, pp. 192–203.
https://doi.org/10.1016/j.actaastro.2016.12.021

[16]. Ho, K., Wang, H., DeTrempe, P. A., Sarton du Jonchay, T., and Tomita, K., “Semi-Analytical Model for Design and Analysis of On-Orbit Servicing Architecture,” *Journal of Spacecraft and Rockets*, Vol. 57, No. 6, 2020, pp. 1129-1138.
https://doi.org/10.2514/1.A34663

[17]. Sears, P., and Ho, K., “Impact Evaluation of In-Space Additive Manufacturing and Recycling Technologies for On-Orbit Servicing,” *Journal of Spacecraft and Rockets*, Vol. 56, No. 6, 2018, pp. 1498–1508.
https://doi.org/10.2514/1.A34135

[18]. Sarton du Jonchay, T., “Modeling and Simulation of Permanent On-Orbit Servicing Infrastructures Dedicated to Modularized Earth-Orbiting Platforms,” *Master’s Thesis, University of Illinois at Urbana-Champaign*, 2017.

[19]. Matos de Carvalho, T. H., and Kingston, J., “Establishing a Framework to Explore the Servicer-Client Relationship in On-orbit Servicing,” *Acta Astronautica*, Vol. 153, 2018, pp. 109–121.
https://doi.org/10.1016/j.actaastro.2018.10.040

[20]. Galabova, K. K., de Weck, O. L., “Economic case for the retirement of geosynchronous communication satellites via space tugs,” *Acta Astronautica*, Vol. 58, No. 9, 2006, pp. 458-498.
https://doi.org/10.1016/j.actaastro.2005.12.014

[21]. Hudson, J. S., Kolosa, D., “Versatile On-Orbit Servicing Mission Design in Geosynchronous Earth Orbit,” *Journal of Spacecraft and Rockets*, Vol. 57, No. 4, 2020, pp. 844-850.
https://doi.org/10.2514/1.A34701

[22]. Sarton du Jonchay, T., Chen, H., Gunasekara, O., and Ho, K., “Framework for Modeling and Optimization of On-Orbit Servicing Operations under Demand Uncertainties,” *Journal of Spacecraft and Rockets*, 2021 (Articles in advance). https://doi.org/10.2514/1.A34978

[23]. Jagannatha, B., and Ho, K., “Event-Driven Network Model for Space Mission Optimization with High-Thrust and Low-Thrust Spacecraft,” *Journal of Spacecraft and Rockets,* Vol. 57, No. 3, pp. 446-463, 2020.
https://doi.org/10.2514/1.A34628

[24]. Vielma, J. P., Ahmed, S., and Nemhauser, G., “Mixed-Integer Models for Nonseparable Piecewise Linear Optimization: Unifying Framework and Extensions,” *Operations Research*, Vol. 58, No. 2, Oct. 2009, pp. 303-315.
http://dx.doi.org/10.1287/opre.1090.0721

[25]. Chen, H., and Ho, K., “Integrated space logistics mission planning and spacecraft design with mixed-integer nonlinear programming,” *Journal of Spacecraft and Rockets*, Vol. 55, No. 2, 2018, pp. 365–381.
https://doi.org/10.2514/1.A33905

[26]. Ishimatsu, T., de Weck, O.L., Hoffman, J.A., Ohkami, Y., and Shishko, R., “Generalized multicommodity network flow model for the earth-moon–mars logistics system,” *Journal of Spacecraft and Rockets*, Vol. 53, No. 1, 2016, pp. 25–38.
https://doi.org/10.2514/1.A33235

[27]. Ho, K., de Weck, O.L., Hoffman, J.A., and Shishko, R., “Dynamic modeling and optimization for space logistics using time-expanded networks,” *Acta Astronautica*, Vol. 105, No. 2, 2014, pp. 428–443.
https://doi.org/10.1016/j.actaastro.2014.10.026

[28]. Ho, K., de Weck, O.L., Hoffman, J.A., and Shishko, R., “Campaign-level dynamic network modelling for spaceflight logistics for the flexible path concept,” *Acta Astronaut*ica, Vol. 123, 2016, pp. 51–61.
https://doi.org/10.1016/j.actaastro.2016.03.006

[29]. Chen, H., Lee, H., and Ho, K., “Space transportation system and mission planning for regular interplanetary missions,” *Journal of Spacecraft and Rockets*, Vol. 56, No. 1, 2019, pp. 12–20.
https://doi.org/10.2514/1.A34168

[30]. Chen, H., Sarton du Jonchay, T., Hou, L., and Ho, K., “Integrated in-situ resource utilization system design and logistics for Mars exploration,” *Acta Astronautica*, Vol. 170, 2020, pp. 80-92.
https://doi.org/10.1016/j.actaastro.2020.01.031

[31]. Beale, E. M. L., and Tomlin, J. A., "Special Facilities in a General Mathematical Programming System for Non-Convex Problems using Ordered Sets of Variables," *Operational Research*, Vol. 69, Nos. 447-454, 1970, Paper 99.

[32]. Session 6: Analytical Approximations for Low Thrust Maneuver. Massachusetts Institute of Technology. https://ocw.mit.edu/courses/aeronautics-and-astronautics/16-522-space-propulsion-spring-2015/lecture-notes/MIT16_522S15_Lecture6.pdf. [Accessed 1/31/2021]

[33]. In-Space Propulsion Data Sheets – Aerojet Rocketdyne, https://www.rocket.com/sites/default/files/documents/In-Space%20Data%20Sheets%204.8.20.pdf [retrieved 8 July 2021].

[34]. UCS Satellite Database. Union of concerned Scientists, www.ucsusa.org/resources/satellite-database. [Accessed 7/17/2020]

[35]. Sullivan, B. R., and Akin, D. L., "Satellite Servicing Opportunities in Geosynchronous Orbit," *AIAA SPACE 2012 Conference and Exposition*, AIAA Paper 2012-5261, 2012. https://doi.org/10.2514/6.2012-5261

[36]. Premiere for Europe: Jules Verne Refuels the ISS. The European Space Agency. http://www.esa.int/Science_Exploration/Human_and_Robotic_Exploration/ATV/Premiere_for_Europe_Jules_Verne_refuels_the_ISS. [Accessed: 9/28/2020].

[37]. Quilty, C., and Moeller, A., "January 2018 Satellite & Space Monthly Review", Quilty Analytics, Feb. 2018. https://www.quiltyanalytics.com/wp-content/uploads/2018_01-Satellite-Monthly-1.pdf [retrieved 8 July 2021].